\documentclass[11pt]{amsart}
\usepackage{amsmath}
\usepackage{amssymb}
\usepackage[lite]{amsrefs}

\usepackage[all]{xypic}
\usepackage[mathscr]{eucal}

\newtheorem{theorem}{Theorem}[section]
\usepackage[mathscr]{eucal}

\newtheorem{lemma}[theorem]{Lemma}
\newtheorem{fact}[theorem]{Fact}
\newtheorem{cor}[theorem]{Corollary}
\newtheorem{claim}[theorem]{Claim}
\newtheorem{remark}[theorem]{Remark}
\newcommand\vlabel{\label}
\newcommand{\la}{\langle}
\newcommand{\ra}{\rangle}

\newcommand{\CN}{{\mathcal N}}

\newcommand{\CR}{{\mathcal R}}
\newcommand{\CM}{{\mathcal M}}

\newcommand{\CI}{\mathcal{I}}
\newcommand{\CF}{{\mathcal F}}
\newcommand{\CG}{{\mathcal G}}

\newcommand{\sub}{\subseteq}
\newcommand{\pf}{\proof }

\newcommand{\CU}{\mathcal U}

\newcommand{\lb}{\lfloor}
\newcommand{\rb}{\rfloor}

\renewcommand{\leq}{\leqslant}
\renewcommand{\geq}{\geqslant}
\newcommand{\tG}{\tilde G}

\newcommand{\reals}{\mathbb R}

\newlength{\intlength}

\title[Central extensions]{On central extensions and definably compact
groups in o-minimal structures}

\date{October 27th 2008 }
\author[Hrushovski]{Ehud~Hrushovski}\address{Department of Mathematics, Hebrew U., Jerusalem, Israel}
\email{ehud@math.huji.ac.il}
\author[Peterzil]{Ya'acov~Peterzil}
\address{Department of
Mathematics, University of Haifa, Haifa, Israel}
\email{kobi@math.haifa.ac.il}
\author[Pillay]{Anand Pillay}\address{Department of Pure Mathematics, University of Leeds, Leeds, LS2 9JT,
U.K. }
\email{pillay@maths.Leeds.ac.uk}\thanks{The first author thanks ISF grant 1048-078, and travel funds from
an EPSRC grant. The third author thanks the European Commission (Marie Curie Excellence chair) and  EPSRC
(grant EP/F009712/1). All authors thank the Marie Curie Modnet network for facilitating the collaboration}

\begin{document}
\maketitle

\begin{abstract} We prove several structural results on definably compact groups $G$ in
o-minimal expansions of real closed fields such as (i) $G$ is definably an almost direct product of a
semisimple group
and a commutative group, (ii) $(G,\cdot)$ is elementarily equivalent to $(G/G^{00},\cdot)$. We also prove
results on the internality of finite covers of $G$ in an o-minimal environment, as well as deducing the full compact
domination conjecture for definably compact groups from the semisimple and commutative cases which were already settled.

These results depend on key theorems about  the interpretability of central and
finite extensions of definable groups, in the o-minimal context. These methods and others also yield interpretability results for universal covers of arbitrary definable real Lie groups.  

\end{abstract}

\section{Introduction and preliminaries}
This paper is motivated partly by questions coming out of our  paper
\cite{nip1}, especially whether, for a definably compact group $G$ in
an o-minimal structure (say expanding a real closed field), $G$ and $G/G^{00}$ are elementarily
equivalent as groups. We solve this problem (see Theorem 7.1) and in the process manage to tie up several
loose ends regarding definable groups in o-minimal structures. For now we will just say ``o-minimal
structure" $\CM$ but often there are additional assumptions on $\CM$ such as expanding a real closed field,
or expanding an ordered group, which appear explicitly in the statements. One of the main results, Theorem
6.1, considers a definably connected central extension $\tilde G$ of a semisimple group $G$ by $A$, all
definable in $\CM$ and says  that the exact sequence
$A \to \tilde G \to G$ of groups, is essentially bi-interpretable with the pair $\la G,A \ra$ of groups.
Corollary 6.2 deduces that any such $\tilde G$ (in particular any definably compact group) is elementarily
equivalent, as
a group, to a semialgebraic real Lie group. From this it is not hard to deduce (Corollary 6.4) that a
definably compact definably connected group is definably an almost direct product of a semisimple group and
a commutative group. Corollary 6.5 strengthens this to central extensions of definably compact semisimple
groups. Section 8 contains interpretability and internality results for finite (not necessarily central)
extensions of groups, again definable in an
o-minimal structure $\CM$. In section 9, definable groups which are not necessarily definably connected are
considered, and Corollary 6.4 (elementary equivalence to semialgebraic Lie groups) is generalized (see
Theorem 9.4).
In section 10 we point out how the compact domination conjecture (for definably compact groups in o-minimal expansions of real closed fields) follows from our results, together with earlier work.

The above results rely on the main technical theorem (Theorem \ref{interpret}) about the  interpretability of central
extensions in an o-minimal context, which appears in section 2. The theorem roughly says that under certain  assumptions, a definable central extension $\tG$ of a definable group $G$ can be interpreted in the
 two-sorted structure $\la G, Z(\tG)\ra$ (possibly, after expanding $G$ by definably connected components of some definable sets). In fact we also note that if the base $o$-minimal structure $\CM$ is an expansion of the reals (in which case we sometimes call a group definable in $\CM$, a definable real Lie group), 
 the assumption that $\tG$ is definable can be omitted, obtaining interpretability results for central topological extensions of suitable definable Lie groups (see section \ref{Omitting} and Theorem 2.8). A version for finite extensions appears in section 8.1
 (see Theorem 8.4). Also in section 8.1 a result with a similar flavour is proved for arbitrary connected definable real Lie groups $G$: for example, the universal cover $\pi:\tG \to G$ is interpretable in the two sorted structure consisting of the given $o$-minimal expansion of $\mathbb R$ together with $\la ker(\pi), + \ra$.

Sections 3, 4 and 5 are devoted to checking
that various hypotheses of Theorem 2.1 hold in the cases we are interested in.

\vspace{5mm}
\noindent
Our notation is in on the whole standard. However, as we are concerned with issues of interpretability in
certain reducts, we will mention the relevant notation.

We do not in general keep track of the definable-interpretable distinction. Hence, when we talk about {\em
definable} sets, groups etc.  in a structure $\CN$, we mean definable (with parameters) in $\CN^{eq}$. We
sometimes say (by abuse of notation) that a set $X$ is $\CN$-definable if $X$ is definable in the above
sense (with parameters) in the structure $\CN$.

In general, $\CM$ is an o-minimal structure, with $M$ is its universe and as  a rule our groups $G$,
$\tilde G$, $H$ etc. are all definable in $\CM$ (again with parameters). However, in order to use results
from \cite{pi} about topology of groups we add the extra assumption that {\em in the structure $\CM$, every
such group is definably isomorphic to a group whose universe is a subset of $M^n$}. All main results assume
that $\CM$ has definable Skolem functions so this assumption is obtained for free in that setting. We may
often want to view the group $G$ as a structure in just the group language, i.e. $\la G,\cdot \ra$, in
which case we sometimes write (with maybe some ambiguity) this structure as just $G$.

\vspace{5mm}
\noindent

We now review some earlier results, mainly about definably simple
and definable semisimple groups.

A \emph{definably simple} group is a definable, non-abelian group
with no definable normal subgroup. A \emph{semisimple} group is a
definable group with no infinite definable normal abelian subgroup
(because of DCC, the definability requirement is superfluous).

We summarize the main results which we will be using here:

\begin{fact}\vlabel{facts}Let $\CM$ be an o-minimal structure.
\begin{enumerate}
\item If $G$ is a definably simple group then there is in $\CM$
a definable real closed field $R$ and a real algebraic group $H$
defined over $R_{alg}\sub R$ the subfield of real algebraic numbers,
such that $G$ is definably isomorphic in $\CM$ to $H(R)^0$, the
definably connected component of $H(R)$ (see \cite{pps1}*{4.1} for the
existence of an algebraic group $H$ and \cite{pps3}*{5.1} and its proof, for the fact
that $H$ can be defined over $R_{alg}$).
\item If $G$ is definably simple then it is either bi-interpretable,
over parameters, with a real closed field or with an algebraically
closed field of characteristic zero, \cite{pps2}.
\item If $G$ is definably connected and semisimple then $Z(G)$ is finite and $G/Z(G)$ is definably
    isomorphic in $\CM$ to
the direct product of definably simple groups (see \cite{pps1}*{4.1}).
\item If $G$ is definably compact and definably connected then either
$G$ is abelian or  $G/Z(G)$ is semisimple, \cite{ps1}*{5.4}.
\item If $G$ is definably simple in a sufficiently saturated
structure then it is abstractly simple if and only if $G$ is not
definably compact, \cite{pps3}*{6.3}.
\end{enumerate}
\end{fact}

Note that the bi-interpretability of (2) above is necessarily over
parameters (see \cite{pps2}*{Remark 4.11}). Also, the Skolem functions
assumption can be omitted if $G$ is definable in the real sort of $\CM$.

\section{The main interpretability theorem}
We recall that $\CM$ is an $o$-minimal structure.

Let $G$ be a definable  group. By a {\em definable central extension
of $G$} we mean the following data: definable groups $A,\tilde G$,
definable homomorphisms: $i:A\to \tilde G$, $\pi:\tilde
G\to G$ with
$$1\rightarrow A\xrightarrow{i} \tilde G\xrightarrow{\pi} G\rightarrow 1$$
exact and $i(A)$ central in $\tilde G$. We let $\la A,\tilde G, G,
i,\pi\ra$ denote the three group structures, together with the maps
$i$ and $\pi$.  We say that the
central sequence above is (definably)
 isomorphic to another definable exact sequence
$$1\rightarrow A_1\xrightarrow{i_1} \tilde G_1\xrightarrow{\pi_1} G_1\rightarrow
1$$ if there are (definable)  group isomorphisms
$h_A:A_1\to A$, $h_{\tilde G}:\tilde G_1\to \tilde G$ and
$h_G:G_1\to G$, which commute with the exact sequence maps.

%We say that  $\la A,\tilde G, G,i,\pi\ra$ is interpretable in a
%theory $T$ if the groups and the corresponding maps are
%interpretable in $T$ (don't we mean here ``interpretable in a
%structure'' and not a theory?).

When $G$ is a definable group whose universe is in the real sort of an an o-minimal $\CM$ then it has a
canonical group topology (see \cite{pi}), with respect to which
every $\CM$-definable subset of $G$ has finitely many definably
connected components. Let $\bf G$ be an  arbitrary expansion of the
group $G$ (not necessarily definable in $\CM$). We say that $\bf G$
{\em has property $\rho$} if for every $\la G,\cdot \ra$-definable
$X\sub G^n$, every definably connected component of $X$ (with
respect to the group topology of $G$) is definable in $\bf G$,
possibly over new parameters.
%{\bf note that I added the
%comment about the parameters. maybe you both realized it all the time. I
%kept getting annoyed by the need for them}.

Given the abelian group $A$, we use  $\la {\bf G}, A\ra$ to
denote the two-sorted structure of the two groups, where $G$ is
equipped with its $\bf G$-structure and $A$ with just its group
structure.

 Recall that for a definable group $H$, the commutator subgroup
$[H,H]$ is a countable union of definable sets, which might not be
definable itself. More precisely, if we denote by $[H,H]_n$ the
definable set of all products of $n$ commutators in $H$, then
$[H,H]=\bigcup_{n\in \mathbb N} [H,H]_n.$

\begin{theorem} \vlabel{interpret} Let $\CM$ be an
o-minimal structure and assume $E= 1\rightarrow A\xrightarrow{i}
\tilde G\xrightarrow{\pi} G\rightarrow 1$ is an $\CM$-definable
central extension of $G$, $\bf G$ an  arbitrary expansion of $G$
such that:

\noindent (1) $\bf G$ has property $\rho$.

 \noindent (2) For every $n$, the set $i(A)\cap [\tilde G,\tilde G]_n$ is
finite (we call this {\em property  $Z$}).

\noindent (3) There
exists $r\in \mathbb N$ with  $G=[G,G]_r$.

Then $E$ (that is the structure $\la A,\tilde G, G, i, \pi \ra)$ can be defined
in $\la {\bf G}, A \ra$ over $A$ and $G$. More precisely, there is an exact sequence
$$E'=1\rightarrow
A\rightarrow \tilde G'\rightarrow G\rightarrow 1,$$ definable in
$\la {\bf G}, A \ra$ over an imaginary parameter $\bar c$, such that $E'$ is
definably isomorphic, in the structure $\la A,\tilde G, {\bf G}, i, \pi \ra)$ (note that $G$ is expanded),
to the
sequence $E$, with $h_A$, $h_G$ the identity maps.

The imaginary parameter $\bar c$ names a map from the set of
definably connected components of a $0$-definable set in
the group $\la G,\cdot\ra$, onto a finite subset of $A$. This map
 is $0$-definable in the
structure $\la A,\tilde G,G,i,\pi\ra$.

\end{theorem}

\proof Note that our assumption implies that $\tilde G$ can be
written as the group product of the subgroups $i(A)$ and $[\tilde
G,\tilde G]_r$.  Our
goal is to produce an $\CM$-definable surjective map $j:G^{2r}\times
A \to \tilde G$, such that the pull-back under $j$ of equality in
$\tilde G$ and of the group operation are both definable in
$\la {\bf G}, A \ra$.

For $x,y\in G$, we let $[x,y]=xyx^{-1}y^{-1}$. For $n\geq 0$, we let
$w_n(x_1,y_1,\ldots, x_n,y_{n})$ be the word in the free group given
by the product of the $n$-commutators $[x_1,y_1]\cdots
[x_{n},y_{n}]$. For any group $K$, we let $F_{n,K}:K^{2n}\to K$ be
the associated function which evaluates the word $w_n$ in $K$. The
image of $K^{2n}$ under $F_{n,K}$ is exactly $[K,K]_n$. We also
have, for $\bar h_1 \in K^{2m}, \bar h_2\in K^{2n}$,
$$F_{m,K}(\bar h_1)\cdot  F_{n,K}(\bar h_2)=F_{m+n,K}(\bar h_1,\bar h_2),$$
and if for $\bar h=(x_1,y_1,\ldots, x_n,y_n)$ we let $inv(\bar
h)=(y_n,x_n,y_{n-1},x_{n-1},\ldots, y_1,x_1)$ then
$$F_{n,K}(\bar h)^{-1}=F_{n,K}(inv (\bar h)).$$

We use $\pi$ to denote the map from $\tilde G^{2n}$ to $G^{2n}$
which is induced by $\pi$ in each coordinate.

\begin{claim}\vlabel{factor}
For every $\bar g_1, \bar g_2\in \tilde
G^{2n}$, if $\pi(\bar g_1)=\pi(\bar g_2)$ then $F_{n,\tilde G}(\bar
g_1)=F_{n,\tilde G}(\bar g_2)$.
\end{claim}
\pf This is immediate from the fact that each  coordinate of $\bar
g_1$ differs from the corresponding coordinate of  $\bar g_2$ by a
central  element of $i(A)$, and on tuples from the center of $\tilde
G$, the map $F_{n,\tilde G}=e$.\qed

It follows from \ref{factor} that there is an $\CM$-definable
surjective map $k_{n}:G^{2n}\to [\tilde G,\tilde G]_n$ such that
$F_{n,\tilde G}$ factors through $\pi$ and $k_{n}$ (see diagram
below) . Also, for $\bar g_1\in G^{2m}, \bar g_2\in G^{2n}$, we have
$$k_m(\bar g_1)\cdot k_n(\bar g_2)=k_{m+n}(\bar g_1,\bar g_2) \mbox{
and } k_n(\bar g_1)^{-1}=k_n(inv(\bar g)).$$
$$
\xymatrix{
\tilde G^{2n} \ar[d]_{\pi} \ar[rr]^{F_{n,\tilde G}} &  &\tilde G \ar[d]^{\pi} \\
G^{2n} \ar[rru]^{k_n} \ar[rr]_{F_{n,G}}&  & G  }
$$

\begin{fact} \vlabel{fact-interpret} (i) The function $k_n:G^{2n}\to \tilde G$ is continuous (here and
below
we always refer to the group topology).

(ii) For every $\bar g_1,\bar g_2\in G^{2n}$, the following are
equivalent:
\begin{enumerate}
\item $F_{n,G}(\bar g_1)=F_{n,G}(\bar g_2)$.
\item $k_n(\bar g_1)$ and $k_n(\bar g_2)$ are in the same
$i(A)$-coset in $\tilde G$.
\item $F_{2n, G}(\bar g_1,inv(\bar g_2))=e.$
\end{enumerate}

(iii) For every $\bar g\in G^{2n}$, $F_{n,G}(\bar g)=e$ if and only if
$k_n(\bar g)\in i(A)$.

\end{fact}
\proof (i) Because the group topology on $G$ equals the quotient
topology (see \cite{b1}*{Theorem 4.3}) and because $F_{n,\tilde G}$ is continuous,
the map $k_{n}$ is also continuous.

(ii) $(1) \Leftrightarrow (2)$: Because $\pi$ is a homomorphism, we
have $\pi k_n(\bar g)=F_{n,G}(\bar g)$, for every $\bar g\in \tilde
G^{2n}$. Hence, $F_{n,G}(\bar g_1)=F_{n,G}(\bar g_2)$ if and only if
$\pi k_n(\bar g_1)=\pi k_n(\bar g_2)$ if and only if $k_n(\bar g_1)$
and $k_n(\bar g_2)$ are in the same $i(A)$-coset.

$(2) \Leftrightarrow (3)$:  $k_n(\bar g_1)\cdot k_n(\bar
g_2)^{-1}\in i(A)$ if and only if $k_{2n}(\bar g_1,inv(\bar g_2))\in
i(A)$ if and only if $F_{2n,G}(\bar g_1,inv(\bar g_2))=e$.

(iii) is immediate from (ii), by taking $\bar g_2=(e,\ldots,e)$.\qed

%\noindent{\bf Claim }
% For every $\bar g_1, \bar g_2\in G^{2n}$,
%$F_{n,G}(\bar g_1)=F_{n,G}(\bar g_2)$ (in $G$) if and only if
%$k_n(\bar g_1)$ and $k_n(\bar g_2)$ differ by an element of
%$i(A)\cap [\tilde G,\tilde G]_{2n}$.

%\proof If $F_{n,G}(\bar g_1)=F_{n,G}(\bar g_2)$ then $k_n(\bar g_1)$
%and $k_n(\bar g_2)$ are in $[\tilde G,\tilde G]_n$ and project to
%the same element in $G$ and therefore differ by an element of
%$i(A)$. The difference of any two elements in $[\tilde G,\tilde
%G]_n$ is clearly in $[\tilde G,\tilde G]_{2n}$. The converse is
%immediate as well.\qed

  For $k\geq 0$, we let
$$G(k)=\{\bar g\in G^{2k}:F_{k,G}(\bar g)=e\}.$$

By our assumption, the set $i(A)\cap[\tilde G,\tilde G]_n$ is finite
for every $n$. Because $k_{n}$ is continuous and surjective on
$[\tilde G,\tilde G]_n$ (since $F_{n,\tilde G}$ is) we have:
\begin{fact}\vlabel{fact2-interpret} The set $k_{n}(G(n))$ equals $i(A)\cap [\tilde G,\tilde G]_n$,
and the function $k_n$ is constant on every definably connected
component of $G(n)$.\end{fact}

Given $n$, let $b_1,\ldots, b_{\ell_{n}}\in A$ be such that
$$\{i(b_1),\ldots, i(b_{\ell_n})\}=i(A)\cap [\tilde G,\tilde
G]_n. $$ We have a corresponding partition of $G(n)$ into relatively
clopen sets $W_{n}(b_1), \ldots, W_n(b_{\ell_n})$, with
$k_n(W_n(b_j))=\{i(b_j)\}$. Each $W_{n}(b_j)$ is a finite union of
components of the set $G(n)$ which is itself $\la
G,\cdot\ra$-definable. Hence, by property $\rho$, each such
$W_n(b_j)$ is  $\bf G$-definable, possibly over some parameters.
\\

\noindent{\bf The interpretation}

 We fix $r$ such that
  $\tilde G=i(A)\cdot [\tilde G,\tilde G]_r$.
\\

%  Because $G$ is definably compact, it follows (see \cite{ }) that $H=G/Z(G)$
%has no infinite normal abelian subgroup, namely it is semisimple.
% Therefore, $H$ is
%definably isomorphic to an almost direct product of definably
%compact semialgebraic, definably simple groups defined over the real
%numbers. Each such group $S$ is elementarily equivalent to a compact
%simple real Lie groups. Hence, there is a $k$ such that every
%element of $S$ is a product of $k$ commutators from $S$. It follows
%that there is an $r$ such that $H=[H,H]_r$ and therefore
%$$G=A\star [G,G]_r.$$ We fix such an $r$.

\noindent {\em The Universe: } Consider the map $j:A\times G^{2r}\to
\tilde G$ defined by $$j(a,\bar g)=i(a)\cdot k_r(\bar g).$$

Because $k_r$ is surjective on $[\tilde G,\tilde G]_r$,  the map $j$
is surjective on $\tilde G$, and we have $j(a_1,\bar g_1)=j(a_2,\bar
g_2)$ if and only if $k_r(\bar g_1)\cdot k_r(\bar
g_2)^{-1}=i(a_1)^{-1}\cdot i(a_2)$ if and only if $$k_{2r}(\bar g_1,
inv(\bar g_2))=i(a_1^{-1}\cdot a_2).$$

Let $$B_{2r}=\{b_1,\ldots,b_{\ell_{2r}}\}\sub A$$ be such that
$i(B_{2r})=i(A)\cap [\tilde G,\tilde G]_{2r},$ and let
$${\mathcal W}_{2r}= \{ W_{2r}(b_1),\ldots,
W_{2r}(b_{\ell_{2r}})\}$$ be the corresponding partition of $G(2r)$,
as given above (namely, $k_{2r}(W(b_j))=i(b_j)$). Let $$c_{2r}:
B_{2r}\to \mathcal W_{2r}$$ be the bijection which sends $b_j$ to
$W(b_j)$. Note that the map $c_{2r}$ is $0$-definable in the
structure $\la A,\tilde G,G,i,\pi\ra$ because $W(2r)$, $B_{2r}$ and
$k_{nr}$ are $0$-definable there.

Consider now the equivalence relation $\sim$ induced on $A\times
G^{2r}$ by the map $j$. It is defined by
$$j(a_1,\bar g_1)=j(a_2,\bar g_2) \Leftrightarrow
 (\bar g_1,inv(\bar g_2))\in c_{2r}(a_1^{-1}a_2).$$

We let $\mathcal U= A\times G^{2r}/\sim$. Notice that the
equivalence relation is $0$-definable in the pure group language of
$A$, the pure group language of $G$, together with a function symbol
for $c_{2r}:B_{2r}\to {\mathcal W}_{2r}$, which in particular names
the finite set $B_{2r}\sub A$. By property $\rho$, the function
$c_{2r}$ itself is definable, over parameters, in  $\la {\bf G},A\ra$
(one way to obtain $c_{2r}$ is by naming each element of $B_{2r}$
and then naming an element in each $W(b_j)\in {\mathcal W_{2r}}$).

We denote by $\lb(a,\bar g)\rb $ the $\sim$-class of $(a,\bar g)$.
\\

\noindent {\em The group operation: } We now consider the pull-back
on $\mathcal U$, via the map $j$, of the group operation from
$\tilde G$: We get (because $i(A)$ is central) for every $(b,\bar
h), (a_1,\bar g_1), (a_2,\bar g_2)\in A\times G^{2r}$,
$$\lb (b,\bar h)\rb=\lb (a_1,\bar g_1)\rb \cdot \lb(a_2,\bar g_2)\rb \Leftrightarrow
b\cdot k_{r}(\bar h)= a_1\cdot a_2\cdot k_{r}(\bar g_1)\cdot
k_{r}(\bar g_2)$$ if and only if $$(\bar g_1,\bar g_2,inv(\bar
h))\in c_{3r}(a_2^{-1}a_1^{-1}).$$

As before, this last expression can be defined using the pure group
structure of $A$ and $G$, and a function symbol for
$c_{3r}:B_{3r}\to {\mathcal W}_{3r}$. The map $c_{3r}$ itself is
defined, over parameters, in $\la {\bf G}, A \ra$.

We thus proved the interpretation of the group $\tilde G$ in
$\la {\bf
G}, A \ra$, over the imaginary parameter $\bar c=(c_{2r},c_{3r})$.
The map $i:A\to \tilde G$ is interpreted by $i(a)=\lb
(a,(e,\cdots,e))\rb$ and the map $\pi:\tilde G\to G$ is interpreted
via $\pi(\lb (a,\bar g)\rb)= F_{r,G}(\bar g).$

We therefore obtain in $\la {\bf G}, A \ra$ a central extension of $G$
which is isomorphic to the original
one, as required. \qed
\\

\begin{cor}\vlabel{corinterpret} Consider the assumptions of Theorem \ref{interpret}
and assume further that $\bf G$ is just the group structure  $G$.
Then $\la A, \tilde G,G,i,\pi\ra$ and $\la G, A \ra$ are
bi-interpretable over parameters. The parameters are  in $G$ and $A$
and they are $0$-definable in $\la A, \tilde G,G,i,\pi\ra$.
\end{cor}
\proof This is immediate from Theorem \ref{interpret} (using the
fact that the isomorphism between the two exact sequences is the
identity when restricted to $A$ and to $G$, in the notation of that
theorem).\qed

\begin{remark}\vlabel{remark} Let us return to the
the imaginary parameters $c_{2r},c_{3r}$ used in Theorem
\ref{interpret}:  The maps $c_{2r}, c_{3r}$ define correspondences
between finite subsets $B_{2r}, B_{3r}\sub A$, respectively, and
definably connected components of
 $0$-definable sets in $\la G,\cdot\ra$. However, in an
o-minimal structure the definably connected components of a
$0$-definable set are themselves $0$-definable and
therefore
 $c_{2r}$ and $c_{3r}$ are $0$-definable
in the structure $$\la M,<,\la G,\cdot\ra, \la A,\cdot, B_{2r},
B_{3r}\ra \ra$$ (by that we mean that we add predicates for $B_{2r}$
and $B_{3r}$). In the special case that $A$ itself is finite the
sets $B_{2r}$ and $B_{3r}$ are $0$-definable in $\la M,<,\la
A,\cdot\ra\ra$, so these predicates can be omitted.
 We will later make
use of this fact.
\end{remark}

%A similar situation will occur in .... {\bf Udi, please add here whatever you had
%in mind. Note that because of the way I structured the paper your
%suggestion, to include something about settings in which we have "words
%with constant differential" does not fit very well}
%\end{remark}

%{\bf I am not sure where to fit Udi's old remark 2.4
%but also i don't understand it.
%It seems to say (please correct me):
 %If $G=[G,G]_k$ and
 %if $E: A\to\tilde G\to_{\pi} G$   is a definable (anywhere) central
 %extension of $G$ which splits abstractly then it splits
% definably.

% I cannot make sense of the proof Udi gave. As far as I can see the formula
% $\phi(\tilde x,x)\sub \tilde G\times G$ which he gave just defines
% the graph of $\pi$ (and not of a section). I am probably missing something.}

Finally let us mention an easy general result on definable splitting, which we do not really use, but
nevertheless
is in the spirit of the other results.
\begin{fact} Suppose that $E : A \to \tilde G \to G$ is a central extension (of abstract) groups, and that
$G=[G,G]_k$  for some $k$. Suppose that $E$ splits abstractly, then it splits definably (in the structure
$\la A,\tilde G,G,i,\pi\ra$).
\end{fact}
\proof By the splitting assumption, $\tilde G$ can be written abstractly as a direct product of $i(A)$ and
a subgroup $H\sub \tilde G$, with $\pi:H\to G$ an isomorphism. It follows that  $[\tilde G,\tilde G]\sub H$
and because $[\tilde G,\tilde G]_k$ projects onto $G$ we have $[\tilde G,\tilde G]_k=[G,G]=H$. In
particular $H=[\tilde G,\tilde G]$ is definable
hence $\tilde G$ split definably.\qed

\subsection{The real case}
\vlabel{Omitting}

 There was actually not much use of o-minimality in
the proof of Theorem \ref{interpret}. Mainly, it was used in order to obtain the canonical
partition of $G(w_{2n})$ into finitely many definably connected components,
on each of which the map $k_{2n}$ is constant. Because of o-minimality this partition
could be read-off just using $G$ (and the definably connected components of
$G$-definable sets), independently of $\tilde G$ and $\pi$. In particular, if we work
over the reals then this assumptions can be partially omitted:

We say that  $$E: 1\rightarrow A\xrightarrow{i}
\tilde G\xrightarrow{\pi} G\rightarrow 1$$ is a topological central
extension of $G$, if $A,\tilde G $ and $G$ are topological groups, and the maps $i$ and $\pi$ are
homomorphisms of topological groups. When $G$ is definable in an o-minimal structure (but possibly not $\tilde G$ and
$A$), we always consider $G$ with its o-minimal topology.

\begin{theorem} \vlabel{real-interpret} Let $\CM$ be an
o-minimal structure over the real numbers, $G$ a definable group
 in $\CM$.  Let $E= 1\rightarrow A\xrightarrow{i}
\tilde G\xrightarrow{\pi} G\rightarrow 1$ be a topological central
extension of $G$ (so $A$, $\tilde G$, $\pi$ and $i$ are not assumed to
be definable), $\bf G$ an  arbitrary expansion of $G$
such that:

\noindent (1) $\bf G$ has property $\rho$.

 \noindent (2) For every $n$, the set $i(A)\cap [\tilde G,\tilde G]_n$ is
finite.

\noindent (3) There
exists $r\in \mathbb N$ with  $G=[G,G]_r$.

Then  there is an exact sequence
$E':1\rightarrow
A\rightarrow \tilde G'\rightarrow G\rightarrow 1,$ definable in the structure
$\la {\bf G}, A \ra$ (over parameters) such that $E'$ is isomorphic (in the group language)  to the
sequence $E$, with $h_A$, $h_G$ the identity maps.

If, moreover, we assume that $A$ is definable in $\CM$ (so the group $\tG'$ is endowed with the o-minimal
topology) and $\CM$  then the isomorphism between $\tG $ and $\tG'$ is also a topological one.

\end{theorem}

Later on, in section 8, we will show how, for finite central extensions, assumption (3) can be omitted. For
now, let's note that the above already implies that every finite topological cover of $SL(2,\reals)$ is
topologically isomorphic to a semialgebraic cover.

\noindent{\em Proof of the theorem} All one has to do is go back to the proof of Theorem \ref{interpret} and see
where definability and o-minimality was used. While the existence of the map $k_n:G^{2n}\to \tilde G$
is just a group theoretic fact (of course now $k_n$ is not definable in $M$), something should
be said about the continuity of $k_n$, mentioned in \ref{fact-interpret}. Indeed, this just follows
from the fact that the continuous homomorphism  $\pi: \tilde G\to G$ is a quotient map. The rest of
\ref{fact-interpret} is just group theoretic.

Proceeding to \ref{fact2-interpret}, we still have $k_n(G(n))=i(A)\cap [\tilde G,\tilde G]_n$ and because
of our assumption, this set is finite, which implies by continuity, that $k_n$ is locally
constant. We now use the fact that over the reals the definably connected components are just
the usual connected components and conclude that $k_n$ is constant on every (definably) connected
component of the ${\bf G}$-definable set $G(n)$. Since the rest of the argument takes place fully in $\la
{\bf G},A\ra$, we obtain in this last structure a definable central exact sequence
$$E: 1\rightarrow A\xrightarrow{i'} \tilde G'\xrightarrow{\pi'} G\rightarrow 1,$$ together with a group
isomorphism $h_{\tG}:\tG'\to \tG$, such that all maps commute (with the identity maps on $A$ and $G$).

Assume now that $A$ itself is definable in the o-minimal structure $\CM$. Let's see why $h_{\tG}$ is a
topological homeomorphism as well.

The group $\tG'$ is obtained as a quotient of $A\times G^{2r}$ by an $\CM$-definable equivalence relation
$\sim$ which is itself the pre-image of equality under the continuous map $j$. The isomorphism $h_{\tG} :\tG'\to G$ is
just the map induced by $j$. Note that by \cite{e1}, the structure
$\CM$ has definable choice functions for subsets of $A\times G$, hence there exists a definable set of representatives $X\sub A\times
G^{2r}$ and a definable bijection $\alpha:\tilde G'\to X$.
By the definition of the o-minimal topology on $\tG'$, the map $\alpha$ is
  continuous on some open subset $U\sub \tG'$ and therefore the composition $j\circ \alpha$, which is just $h_{\tG}$, is continuous on $U$ as well. Since $h_{\tG}$ is a group isomorphism it must be continuous everywhere.  Because $\tG'$ is locally compact (and $\tG$ is Hausdorff) the inverse map is continuous as well. \qed
\\

%By the assumption of abstract splitting $E$ is isomorphic to $$E':  A \to A\times G \to G$$  over $A,G$  %(where
%$E'$ has the natural maps  $a \to (a,1)$, $(a,g) \to g$).
%Regarding $E$ as a structure, consider the formula
%$\phi(x,y)$ which says that there exist $y_{1},...,y_{2k}\in \tilde G$ and $x_{1},..,x_{2k}\in G$,
%such that $x = [x_{1},x_{2}]....[x_{2k-1},x_{2k}]$, $y = [y_{1},y_{2}]....[y_{2k-1},y_{2k}]$ and
%$\pi(y_{i}) = x_{i}$ for $i=1,..,2k$.
%Interpreted in $E'$ (with $\tilde G$ replaced by $A\times G$) $\phi(x,y)$ clearly defines the map taking %$g\in G$ to $(1,g)\in A\times G$. As $E'$ is isomorphic to $E$, the formula must also define a splitting of % $E$.

 \noindent
 In the next section, we investigate each of the three assumptions of Theorem
\ref{interpret}.

\section{Perfect groups}
In this section $\CM$ can be taken to be an arbitrary o-minimal
structure. Recall that $G$ is perfect if $[G,G]=G$.
\begin{claim}\vlabel{perfect1}Let $G$ be a definable group.

\noindent (i) If $G$ is perfect then every homomorphic image of $G$
is perfect.

\noindent (ii) The direct product of perfect groups is perfect.

\noindent (iii) (Assume that $\la G,\cdot\ra$ is sufficiently
saturated). If $\tilde G$ is definably connected, $G$ perfect and $\pi:\tilde
G\to G$ is a finite extension then $\tilde G$ is perfect as well.

\noindent (iv) If $G$ is definably simple then it is perfect.

\noindent (v) If $G$ is semisimple and definably connected then it is perfect.

\end{claim}
\proof (i) and (ii) are easy. For (iii), the assumption implies that
$\pi [\tilde G,\tilde G]=G$, and hence $\tilde G=F\cdot [\tilde
G,\tilde G]$ for some finite group $F\sub \tilde G$. It follows that
 $[\tilde G,\tilde G]$ is a $\bigvee$-definable group (i.e. a countable union of definable sets)
 of finite
index, hence its complement is also $\bigvee$-definable. This
implies, using saturation, that $[\tilde G,\tilde G]$ is a definable
subgroup of finite index, contradicting connectedness.

(iv) If $G$ is definably compact and definably simple then $G$ is
elementarily equivalent to a compact simple real Lie group $H$,
by \ref{facts}(1). By topological compactness, there exists an $r$ such
that $[H,H]_r=H$. This is now true for $G$ as well.

 If $G$ is not definably
compact then by  \ref{facts}(5) it is abstractly simple and
therefore $[G,G]=G$.

(v) Assume that $G$ is semisimple and definably connected. Then
$G/Z(G)$ is centerless, definably connected and semisimple. By
\ref{facts}(3), the group $G/Z(G)$ can be written as a direct
product of $\CM$-definable definably simple groups. the result now
follows from (iv), (ii) and (iii). \qed

\section{ Property $\rho$ and semisimple groups}

In this section we assume that $\CM$ is an arbitrary o-minimal
structure. Recall that an expansion ${\bf G}$ of an $\CM$-definable
group $G$ is said to have property $\rho$ if the definably connected
components of every $\la G,\cdot\ra$-definable subset of $G^n$ are
definable in ${\bf G}$ (possibly over new parameters).

\begin{claim}\vlabel{T1.2} Assume that $G_1,\ldots, G_k$ are definable groups,
such that the theory of each pure group $\la G_i,\cdot\ra$ satisfies
property $\rho$. Then the theory of the pure group
$G=G_1\times\cdots \times G_k$, expanded by a predicate for every
$G_i\sub G$, satisfies property $\rho$ as well.
\end{claim}
\proof Every $\la G,\cdot\ra$-definable set $X\sub G^n$ is a finite
union of sets of the form $X_1\times \cdots \times X_k$, where $X_i$
is a $G_i$-definable subset of $G_i^n$. By our assumptions on each
$G_i$, we may assume that each $X_i$ is definably connected. Each
definably connected component of $X$ is a finite union of such
cartesian products and therefore $\la G,\cdot\ra$-definable, together with
predicates for every $G_i$.\qed

\begin{lemma} \vlabel{T1.1} If $G$ is a definably simple
 group then the pure group $\la G, \cdot\ra$ has property
$\rho$.
\end{lemma}
\proof By \ref{facts}(2), there are two cases to consider: If $G$ is
unstable then it is a semialgebraic group which is bi-interpretable
(over parameters) with a real closed field. It follows that every
$\la G,\cdot\ra$-definable set $X\sub G^n$ is semialgebraic and
every definably connected component of $X$ is again semialgebraic,
and therefore $\la G,\cdot\ra$-definable (possibly over parameters).

  If $G$ is stable then it
is a linear algebraic group over a definable algebraically closed
field $K$. Because $K$ is a definable  algebraically closed field in
the o-minimal structure $\CM$, then, by \cite{ps}, a maximal real
closed subfield $R\sub K$ is definable in $\CM$ and we have
$K=R(\sqrt{-1})$. Since $G$ is a linear algebraic group over $K$, we
may assume that $G\sub K^{\ell}$ for some $\ell$ and that its
group-topology agrees with that of $K^{\ell}$, identified with
$R^{2\ell}$. In particular, the definably connected components of
every definable subset of $G^n$ in the sense of the group topology
are the same as those in the sense of the field $R$.

By \ref{facts}, $G$  is bi-interpretable (again over parameters)
with $K$ and hence the $\la G,\cdot\ra$-definable subsets of $G^n$
are exactly the $K$-constructible sets. It is therefore sufficient
to prove:
\begin{claim}\vlabel{construct} If $K=R(\sqrt{-1})$ is an algebraically closed field definable in
an o-minimal $\CM$ and $X\sub K^n$ is a $K$-constructible set then
every definably connected component of $X$ (in the sense of $R$) is
$K$-constructible.
\end{claim}
\proof The set $X$ is of the form $$X=\bigcup_{i=1}^r(X_i\setminus
Y_i),$$ with each $X_i$ an irreducible algebraic variety and
$Y_i\sub X_i$ an algebraic variety of smaller algebraic dimension.
We claim that each $X_i\setminus Y_i$ is definably connected.

Indeed, it is known that if $V\sub {\mathbb C}^n$ is an irreducible
complex variety then $Reg(V)$ (the set of complex regular points of
$V$) is a connected set, dense in $V$. If we now work in the
structure $\la \reals, <,+,\cdot, \mathbb \ra$ then, by quantifying
over parameters, this fact carries over to $\la R,<,+,\cdot,K\ra$.
Thus, every $Reg(X_i)$ is definably connected in the sense of $R$
and dense in $X_i$.

Thus, the set $Reg(X_i)$ is a definably connected $R$-manifold of
even $R$-dimension $2k$, and we have $\dim_R(Y_i)\leq 2k-2$ (we let
$\dim_R(Y_i)$ denote the o-minimal dimension of $Y_i$ with respect
to $R$, which is twice the algebraic dimension of $Y_i$). The set
$Reg(X_i)\setminus Y_i$ is therefore still definably connected,
dense in $Reg(X_i)$ and so, also dense in $X_i\setminus Y_i$. It
follows that $X_i\setminus Y_i$ is definably connected.

Finally, each definably connected component of $X$ must be a finite
union of sets of the form $X_i\setminus Y_i$, so constructible.

With this ends the proof of the claim and of Lemma \ref{T1.1}\qed

\vspace{5mm}
\noindent
Part (ii) of the theorem below follows from a general result by
Edmundo, Jones and Peatfield, see \cite{ejp}.

\begin{theorem}\vlabel{T1.3} If $\tilde G$ is semisimple and definably connected
then
 \\(i) $\la \tilde G,\cdot\ra $ is bi-interpretable with  $\la \tilde G/Z(\tilde G),\cdot\ra $,
  after naming a parameter $\bar b$ from $\tilde G/Z(\tilde G)$. The parameter can be chosen in
$dcl_{\CN}(\emptyset)$, for $\CN=\la M,<,\la G,\cdot\ra \ra$ and also in $dcl_{\tilde G}(\emptyset)$, where
$\tilde G=\la \tilde G,\cdot\ra$.
\\(ii) There is an $\CM$-definable real closed field $R$, such that
the  group $\tilde G$ is definably isomorphic in $\CM$, over
parameters, to a semialgebraic group $\tilde G'$ over the field of
real algebraic numbers $R_{alg}\sub R$.
\\(iii) $\tilde G$ has property $\rho$.
\end{theorem}
\proof Consider the extension
$$Z(\tilde G)\xrightarrow{\i}\tilde G\xrightarrow{\pi} G.$$
\\

\noindent {\bf Case I} $G$ is definably simple.

(i) By \ref{T1.1}, the pure group $G$ has property $\rho$ and by
\ref{perfect1}, it is a perfect group. Clearly $\tilde G$ has
finite center and hence has property $Z$.  Therefore,  by Corollary
\ref{corinterpret}, the structure $\la Z(\tilde G),\tilde
G,G,i,\pi\ra$ is bi-interpretable with $\la \la G, \cdot\ra, \la
Z(\tilde G),\cdot\ra\ra $, after naming finitely many elements in $G$
and in $Z(\tilde G)$. Moreover, the finite group $Z(\tilde G)$
itself can be defined in $G$ after naming finitely many elements
(which we may even assume to belong to some $0$-definable finite
subgroup of $G$). Thus, $\la Z(\tilde G),\tilde G,G,i,\pi\ra$ is
bi-interpretable, over parameters in $G$, with $\la G,\cdot\ra$.  Finally, we note that $\la Z(\tilde
G),\tilde
G,G,i,\pi\ra$ is bi-interpretable with $\la \tilde G,\cdot\ra$,
hence $\la \tilde G,\cdot\ra$ is bi-interpretable over parameters,
with $\la G,\cdot\ra$. As
was observed in Remark \ref{remark}, the parameters which we use can
be taken in $dcl_{\CN}(\emptyset)$, where $\CN=\la M,<,\la
G,\cdot\ra \ra$.

(ii) By \ref{facts}(1), we may assume that $G$ is an
$R$-semialgebraic group, defined over $R_{alg}\sub R$, for some
$\CM$-definable real closed field $R$. As was shown above, $\la
\tilde G,\cdot\ra$ is interpretable, now over the empty set, in $\la
R,<,\la G,\cdot\ra \ra$. In particular, the group $\tilde G$ is
definably isomorphic in $\CM$ to a group $H$ which is interpretable
in the field $R$, over $R_{alg}$. By elimination of imaginaries in
real closed fields, $H$ can be chosen to be definable.

(iii) We now want to show that $\tilde G$ has property $\rho$. If $G$ is
stable then, by \ref{facts}(2), it is bi-interpretable, over
parameters, with an algebraically closed field $K$. In this case,
because $\tilde G$ and $G$ are bi-interpretable, $\tilde G$ is
definably isomorphic, over parameters, to an algebraic group $H$
over $K$. By Claim \ref{construct}, if $X\sub H^n$ is constructible
over $K$  then its definably connected components are constructible
over $K$ as well. Because of the bi-interpretability of $G$ (so also
of $H$) with $K$, these components are definable in $\la
H,\cdot\ra$, so $H$ (hence $\tilde G$) has property $\rho$.

If $G$ is unstable then it is bi-interpretable with a real closed
field $R$ and therefore, by (i), $\tilde G$ is also bi-interpretable with $R$.
 By (ii), we may assume that $\tilde G$ is semialgebraic over $R_{alg}$.
 This implies that every semialgebraic subset of $\tilde G^n$ is definable
 in the pure group $\la \tilde G,\cdot\ra$. In particular, $\tilde G$ has property $\rho$.
\\

\noindent{\bf Case II} $G$ is semisimple.

(i) We first claim that $\tilde G$ is bi-interpretable with
$G=\tilde G/Z(\tilde G)$ together, possibly with finitely many
constants. For that, we need to establish the three assumptions of
Theorem \ref{interpret} (with $A=Z(\tilde G)$):

By \ref{facts}, $G$ is definably isomorphic in $\CM$ to a product
$H_1\times\cdots\times H_k$, of definably simple groups. Each of the
$H_i'$s is the centralizer of the other $k-1$ groups hence it is
definable in the pure group langauge of $G$ (after naming
parameters). It follows from \ref{T1.1} and \ref{T1.2} that $G$ has
property $\rho$. By \ref{perfect1}, $\tilde G$ is perfect. Because
$Z(\tilde G)$ is finite, we clearly have property $Z$. We can now
apply Theorem \ref{interpret} exactly as in Case I.

(ii) This is identical to the proof in Case I.

(iii) We claim that $\tilde G$ has property $\rho$:

Assume that $G=H_1\times\cdots\times H_k$, for definably simple
$H_i$'s and let $\tilde H_i$ be the pull back of $H_i$ under the
inverse image of $\pi:\tilde G\to G$. Each $\tilde H_i$ is a finite
central extension of $H_i$, and if we let $\tilde H=\tilde H_1\times
\cdots\times \tilde H_k$ and $\tilde \pi:\tilde H\to G$ be the
natural projection, then $\tilde \pi$ factors through the finite
extensions $\pi' :\tilde H\to \tilde G$ and $\pi:\tilde G\to G$.

%
%(11)
$$
\xymatrix{
\tilde H \ar[rd]_{\tilde \pi} \ar[rr]^{\pi'} & & \tilde G \ar[ld]^{\pi} \\
& G &
}
$$

\noindent {\bf Claim} $\tilde H$ has property $\rho$.

\proof  Each $\tilde H_i$ is a finite central extension of a definably simple group
so by Case I, it has property
$\rho$.  Therefore, by \ref{T1.2}, in order to see
that $\tilde H$  itself has property $\rho$ it is sufficient to see
that each $\tilde H_i$, $i=1,\ldots,k$, is definable in the pure
group $\la \tilde H,\cdot\ra$, possibly with parameters. Let us see
why $\tilde H_1$ is definable. First note that the centralizer of
$\tilde H_2\cup\cdots\cup  \tilde H_k$, call it $Z_1$, is $\tilde
H_1\cdot Z(\tilde H)$ (where $Z(\tilde H)$ is finite). The group
 $Z_1$ is a definable, possibly disconnected, group and it is
sufficient to see that we can define in the pure group structure,
its connected component. By \ref{perfect1}(v), $\tilde H_1$ is
perfect, and it is easy to see that $[Z_1,Z_1]\sub \tilde H_1$.
Hence, $\tilde H_1=[Z_1,Z_1]_k$ for some $k$ and this last group is
clearly definable. End of Claim.

 Because $G$ is perfect and has property $\rho$, we can apply
Theorem \ref{interpret} to the finite central extension $\tilde
\pi:\tilde H\to G$ and conclude that the pure groups $G$ and $\tilde
H$ are bi-interpretable, after naming constants from $G$.

Let $X\sub \tilde G^n$ be a $\tilde G$-definable set and let
$X_1,\ldots, X_k$ be its definably connected components, with
respect to the $\tilde G$-topology. Because $\tilde G$ and $ \pi':
\tilde H\to \tilde G$ are  definable in $\la \tilde H, \cdot\ra$ and
continuous, the set $Y=\pi'^{-1}(X)$ is definable in $\tilde H$ and
each $\pi'^{-1}(X_i)$ is a finite union of definably connected
components of $Y$, hence definable in $\tilde H$ (because $\tilde H$
has property $\rho$). It follows that each $X_i$ is definable in
$\tilde H$. However, as we already saw, $\tilde H$ and $G$ are
bi-interpretable and $G$ and $\tilde G$ are bi-interpretable as
well, and therefore each $X_i$ is definable in $\tilde G$ (after
possibly naming finite many parameters). Hence, $\tilde G$ has
property $\rho$.\qed
\\

%{\bf None of you commented on the above proof of (iii). Does it mean that
%we leave the above argument as it is now? I myself could not simplify it}
% Yes (AP)

% {\bf Remark: The last
%argument might seem too complicated. I just could not see a direct
%argument for the following result: Assume that $G_1$ and $G_2$ are
%two definable groups which are bi-interpretable. Then, $G_1$ has the
%$\rho$-property iff $G_2$ has it. The reason is that the group
%topology of $G_1$ is not just the induced topology from $G_2$.}
%\\

\section{Property $Z$}

 Assume now that $\CM$ expands a real closed
field $R$, in a neighborhood of the identity of a definable group $G$.
We denote by $\mathcal G$ its Lie algebra whose
underlying $R$-vector space is the tangent space of $G$ at $e$,
$T_e(G)$ . We recall some facts about groups and Lie algebras, as
presented in \cite{pps1}.

Assume that $G$ is definably connected. To every definable subgroup
$H\sub G$ there is an associated Lie subalgebra ${\frak h}\sub
{\mathcal G}$. The subgroup $H$ is normal in $G$ if and only if
$\frak h$  is an ideal in $\mathcal G$ (see \cite{pps1}*{Theorem 2.32}).
For every $g\in G$, we denote by $Ad_g:T_e(G)\to
T_e(G)$ the differential of the inner automorphism $a_g:x\mapsto
gxg^{-1}$. If ${\mathcal G}_1$ is a linear subspace of $\mathcal G$
then $\mathcal G_1$ is an ideal if and only if it is invariant under
$Ad_g$ for all $g\in G$. (See Claim 2.31 there).

If $H$ is a locally definable subgroup of $G$ (e.g. when $H=[G,G]$)
then, just as for definable subgroups, one can associate to $H$ a
Lie subalgebra $L(H)\sub \CG$.  If $H$ is normal in $G$ then $L(H)$
is an ideal in $\CG$ (indeed, because the whole analysis is local in
nature, it is enough to consider $H$ at a neighborhood of $e$ and in
this case the arguments work as in the definable case).

Recall that a subalgebra $\mathcal A\sub \mathcal G$ is called
central if for every every $\xi\in \mathcal A$ and  $\eta \in
\mathcal G$, we have $[\xi,\eta]=0$. Equivalently (see \cite{pps1}*{Corollary 2.32}), for every $g\in G$,
$Ad_g|\mathcal A=id$.

\begin{theorem}\vlabel{alg1} Let $G$ be a definably connected group. Assume that
 $\mathcal G= \mathcal A+\mathcal \CI$ for a central subalgebra $\mathcal
A$, and an ideal $\mathcal I\sub \CG$. Then $L([G,G])\sub \CI$.
\end{theorem}
\proof We first introduce some notation. Let $\frak{h}=L([G,G])$ and for
$g\in G$, let $\ell_g:G\to G$ be left-multiplication by $g$ and
$r_g:G\to G$ be right-multiplication by $g$. For every $g\in G$ we
have
$$T_g(G)=d_e(\ell_g)(T_e(G))=d_e(r_g)(T_e(G)),$$ and similarly, for
every $h\in [G,G]$ we have
$$T_h([G,G])=d_e(\ell_g)(\frak h)=d_e(r_g)(\frak h).$$

It is therefore sufficient to prove, for some $h\in [G,G]$, that
$d_h(r_h^{-1})(T_h([G,G]))\sub \mathcal I$.

Because $[G,G]$ is locally-definable there exists an $n$ and an open
neighborhood $U\sub G$ of $e$ such that $U\cap [G,G]=U\cap [G,G]_n$.
Consider the function $F_n=F_{n,G}:G^{2n}\to G$ as given earlier, by
the product of $n$-many group commutators. It is not hard to see
that for sufficiently generic $\bar g=(g_1,\ldots, g_{2n})$ in $
F_n^{-1}(U\cap [G,G])$, we have
$$d_{\bar g}(F_n)(T_{g_1}(G)\times \cdots \times T_{g_{2n}}(G))=
T_{F_n(\bar g)}([G,G]).$$

Using the chain rule, it is sufficient to prove that for every $\bar
g\in G^{2n}$, we have $$d_{\bar g}(r_{F_n(\bar g)}^{-1}\circ
F_n)(T_{g_1}(G)\times \cdots \times T_{g_{2n}}(G))\sub \mathcal I.$$

%or equivalently (letting $\bar e=(e_,\ldots, e)\in G^{2n}$ and
%$\ell_{\bar g}=(\ell_{g_1},\ldots,\ell_{g_{2n}}) )
%$$d_{\bar e}(r_{F_n(\bar g)}\circ
%F_n\circ \ell_{\bar g} )(T_{e}(G)^{2n}\sub \mathcal B.$$

We are going to prove this by induction on $n$. For that purpose,
let us a call definable function $\alpha :G^k\to G$ {\em good at
$\bar g\in G^k$} if it satisfies
$$d_{\bar g}(r_{\alpha (\bar g)^{-1}}\circ \alpha)(T_{g_1}(G)\times
\cdots \times T_{g_{k}}(G))\sub \mathcal I.$$

\noindent{\bf Claim} \noindent (i) If $\alpha:G^k\to G$ is good at
$\bar g$ and $\alpha(\bar g)=e$ then for every $h\in G$, the
function  $a_g\circ \alpha$ is good (recall $a_g(x)=gxg^{-1}$).
\\

\noindent (ii) If $\alpha, \beta :G^k\to G$ are good at $\bar g\in
G^k$  then so is $\alpha\cdot \beta$, the group product of the two.
\\

\noindent {\em Proof } (i) This is immediate from the fact that
$\mathcal I$ is invariant under $d_e(a_g)=Ad_g$.

(ii)   If  $\mu:G\times G\to G$ is the group product then
$d_{(e,e)}(\mu)=(id, id)$.  Now, in the special case that
$\alpha(\bar g)=\beta(\bar g)=e$ we have, by the chain rule,
$d_{\bar g}(\mu(\alpha, \beta))=d_{\bar g}(\alpha)+d_{\bar
g}(\beta),$ and therefore $\mu(\alpha,\beta)$ is good at $\bar g$.

More generally, if $\alpha(\bar g)=h_1, \beta(\bar g)=h_2$ then
 $$\alpha(\bar x)
\beta(\bar x)  h_2^{-1}h_1^{-1}=\alpha(\bar x)
h_1^{-1}(h_1(\beta(\bar x) h_2^{-1})h_1^{-1}),$$ and hence
$$r_{(h_1h_2)}^{-1}\circ \mu(\alpha,\beta)=\mu(r_{h_1^{-1}}\circ
\alpha, a_{h_1}\circ  r_{h_2^{-1}}\circ \beta ).$$

By definition, $r_{h_1^{-1}}\circ \alpha$ and  $r_{h_2^{-1}}\circ
\beta$ are good at $\bar g$ and the two functions send $ \bar g$ to
$e$. By (i) and the special case we just did, $\mu(\alpha, \beta)$
is good at $\bar g$ as well. End of Claim.
\\

Because every $F_n$ is a product of commutators, it is sufficient,
using Claim (ii) above, to prove that $F_1(x,y)=xyx^{-1}y^{-1}$ is
good at every $(g,h)\in G^2$. Because $F_1(g,h)=ghg^{-1}h^{-1}$, we
need to show that $r_{hgh^{-1}g^{-1}}\circ F_1$ is good at $(g,h)$,
or equivalently, that $\sigma(x,y)=r_{hgh^{-1}g^{-1}}\circ F_1(gx, hy)$
is good at $(e,e)$. Rewriting $\sigma(x,y)$ we get:
$$\begin{array}{ll}gxhyx^{-1}g^{-1}y^{-1}h^{-1}hgh^{-1}g^{-1}& =gxg^{-1}gh
(yx^{-1}g^{-1}y^{-1}g)(gh)^{-1}\\ & =a_g(x) \cdot a_{gh}(y\cdot
x^{-1}\cdot a_{g^{-1}}(y^{-1}))\mbox{.}\end{array}$$ The right-most
expression can be re-written as
$$a_g(x)\cdot a_{gh}(y)\cdot a_{gh}(x)^{-1}\cdot a_{ghg}^{-1}(y)^{-1}.$$

We have a product of four functions,  each sending $e$ to $e$.
Taking the differential and applying the chain rule we obtain, for
every $u,v\in T_e(G)$:
$$d_{(e,e)}\sigma
(u,v)=Ad_g(u)+Ad_{gh}(v)-Ad_{gh}(u)-Ad_{ghg^{-1}}(v).$$

We now return to our assumptions. Every $u\in \mathcal G$ can be
written as $u=u_1+u_2$, where $u_1\in \mathcal A$ and $u_2\in
\mathcal I$. Because $\mathcal A$ is central we have $Ad_g(u_1)=u_1$
 for every $g\in G$. Hence, $d_{(e,e)}\sigma (u_1+u_2, v_1+v_2)$
equals to:
$$\begin{array}{l}u_1+Ad_g(u_2)+v_1+Ad_{gh}(v_2)-u_1-Ad_{gh}(u_2)-v_1-Ad_{ghg^{-1}}(v_2)=
\\ =Ad_g(u_2)+Ad_{gh}(v_2)- Ad_{gh}(u_2)-Ad_{ghg^{-1}}(v_2)\mbox{.}\end{array}$$
 Because
$\mathcal I$ is invariant under every $Ad_g$ and under $+$, the sum
on the right belongs to $\mathcal I$.\qed

\begin{cor}\vlabel{cor1.1} Let $G$ be a definably connected group, $A\sub  G$
a definable central subgroup and let $\mathcal A$ be the Lie algebra
of $A$. Assume that $\mathcal G$ can be written as a direct sum
$\mathcal G=\mathcal A\oplus \mathcal I$, for some ideal $\CI$. Then
for every $n$, $A\cap [G,G]_n$ is finite.
\end{cor}
\proof The Lie algebra $\mathcal A$ is central in $\mathcal G$ (see \cite{pps1}*{Claim 2.32}). By our
assumption, and Theorem
\ref{alg1} we have $L([G,G])\sub \mathcal I$ and therefore $\mathcal
A\cap L([G,G])=\{0\}$. Because $ A\cap [G,G]$ is a locally definable
group it has a Lie Algebra of the same dimension, which equals
$ \mathcal A\cap L([G,G])$. Hence, $\dim (A\cap[G,G] )=0$ and
therefore $ A\cap[G,G]_n$ is finite for every $n$.\qed

\begin{cor}\vlabel{cor1.2} Let $\tilde G$ be a definably connected
central  extension of a semisimple group $G$, with $L(\tilde G)
=\tilde \CG$. Then

\noindent (i)  For every $n$, the set $Z(\tilde G)\cap [\tilde
G,\tilde G]_n$ is finite.

\noindent (ii) The Lie algebra of the locally definably group
$[\tilde G,\tilde G]$ equals to $[\tilde{\mathcal G},\tilde
{\mathcal G}]$ and we have $\tilde{\mathcal G}=\mathcal Z\oplus
L([\tilde G,\tilde G])$, where $\mathcal Z=L(Z(\tilde G))$.
Moreover,
$$L([\tilde G,\tilde G])\simeq L(G).$$

\end{cor}
\proof (i) First note that $\mathcal Z=L(Z(\tilde {G}))$ is
the center of $\tilde{\mathcal G}$ (see \cite{pps1}*{Claim 2.31})
and the Lie algebra of $\tilde G/Z(\tilde G)$ equals $\tilde {
\mathcal G}/\mathcal Z$. It follows that $\tilde {\mathcal
G}/\mathcal Z$ is a semisimple Lie algebra (see Theorem 2.34
there).

By the Levi decomposition theorem for Lie algebras, $\tilde
{\mathcal G}$ is the semi-direct product of its solvable radical and
a semisimple Lie sub-algebra $\frak h$. Because $\tilde{\mathcal
G}/\mathcal Z$ is semisimple it follows that $\mathcal Z$ is the
solvable radical of $\tilde{\mathcal G}$. We claim that $\frak h=[\tilde{\mathcal G},\tilde{\mathcal G}]$.

Indeed, for $\xi_i=\xi_i+\eta_i \in \tilde{\mathcal G}$, $i=1,2$,
and $\xi_i\in \mathcal Z$ and $\eta_i\in {\frak h}$, we have
$$[\xi_1,\xi_2]=[\xi_1,\xi_2+\eta_2]+[\eta_1,\xi_2]+[\eta_1,\eta_2]=[\eta_1,\eta_2]\in\frak h.$$
It follows that $[\tilde{\mathcal G},\tilde{\mathcal G}]\sub
\frak h$ and because $\frak h$ is semisimple we also have
$\frak h=[\frak h,\frak h]\sub [\tilde{\mathcal
G},\tilde{\mathcal G}]$, and hence $\frak h=[\tilde{\mathcal
G},\tilde{\mathcal G}]$. Therefore, $\frak h$ is an ideal in
$\tilde{\mathcal G}$ and we have $\tilde{\mathcal G}=\mathcal
Z\oplus [\tilde{\mathcal G},\tilde{\mathcal G}]$.

 We can now apply Corollary
\ref{cor1.1} and conclude that $Z(\tilde G)\cap [\tilde G,\tilde
G]_n$ is finite for every $n$.

(ii) By dimension considerations, the above implies that $\dim
L([\tilde{\mathcal G},\tilde{ \mathcal G}])=\dim([\tilde G,\tilde
G])=\dim (\tilde{\mathcal G})-\dim Z(\tilde G)=\dim [\tilde{\mathcal
G},\tilde{\mathcal G}]$. By Theorem \ref{alg1}, $L([\mathcal
G,\mathcal G])\sub [\tilde{\mathcal G},\tilde{\mathcal G}]$, and
therefore $L([\mathcal G,\mathcal G])= [\tilde{\mathcal
G},\tilde{\mathcal G}]$.

Again, by dimension considerations, $\dim ([\tilde G,\tilde G])=\dim G$ and
hence $d\pi$ is an isomorphism of $L([\tilde G,\tilde G])$ and
$L(G)$. \qed

Because every definably compact definably connected group is a
central extension of a semisimple group (\ref{facts}) we
immediately conclude the result below. As we will later see
(Corollary \ref{structure}), this is only a first approximation to
the stronger theorem about the commutator subgroup of a definably
compact group.

\begin{cor} Let $\tilde G$ be a definably connected definably
compact group. Then for every $n$, $Z(\tilde G)\cap [\tilde G,\tilde
G]_n$ is finite.
\end{cor}

\subsection{Omitting the real closed field assumption}

 The real closed field assumption was of course necessary for the
discussion in the last section, because it involved the Lie algebra
of $\tilde G$ at $e$. However, once used this assumption
can be weakened, at least in the definably compact case.

We first recall some notions: An o-minimal expansion of an ordered
group is called {\em semi-bounded} if there is no definable
bijection between bounded and unbounded intervals. There are three
different possibilities for an o-minimal expansion $\CM$ of an
ordered group (see discussion in \cite{pe}):

\noindent 1. $Th(\CM)$ is linear, i.e. $\CM$ is elementarily equivalent to an ordered  reduct of an
ordered
vector space over an ordered division ring.

 \noindent 2. $Th(\CM)$ is not linear and not semi-bounded, in which case
there exists a definable real closed field whose domain is $M$.

\noindent 3. $\CM$ is semi-bounded and $Th(\CM)$ is not linear.

We can now state the following generalization of Corollary
\ref{cor1.2}:
\begin{cor}\vlabel{cor1.4} Let $\CM$ be an o-minimal expansion of an
ordered group and let $G$ be an $\CM$-definable, definably compact,
definably connected group. Then for every $n\in \mathbb N$,
$Z(G)\cap [G,G]_n$ is finite.
\end{cor}
\proof We need to examine the three cases above. If $\CM$ is elementarily equivalent to a reduct
of an ordered vector space then every definable group is
abelian-by-finite (indeed, by \cite{ps1}, if not then a field is
definable in  $\CM$. It is easy to see that this is impossible), so
$G$ is abelian.

If $\CM$ is not linear and not semi-bounded then $\CM$ expands a real
closed field and therefore Corollary \ref{cor1.2} applies. We are
thus left with the semi-bounded nonlinear case. Recall the following from
\cite{pe}:

%{\bf Udi put a "$*$" before and after this paragraph. What was it
%about?}
% I don't know. AP

If $\CM$ is semi-bounded and $Th(\CM)$ is not linear then there exists $\CN\succ \CM$
and an o-minimal expansion $\widehat \CN$ of $\CN$ (by
``expansion'' we mean here that every definable subset of $\CN$ is
definable in $\widehat \CN$, possibly with additional parameters)
and an elementary substructure $\widehat D \prec \widehat \CN$ such
that every interval in $\widehat D$ admits a definable real closed
field. Given $G$,  an $\CM$-definable, definably connected and
definably compact group, we can view $G$ as an $\CN$-definable (and therefore also $\widehat
\CN$-definable) group. Hence, there exists in $\widehat \CN$ a
$0$-definable family $\CF=\{G_s:s\in S\}$ of groups, all
definably compact and definably connected such that $G=G_{s_0}$ for
some $s_0\in S$.  Furthermore, the domain of every such $G_s$ is a
bounded subset of $\widehat \CN$ (see \cite{pe}). Given a fixed
$n\in \mathbb N$, we can now argue in $\widehat D$:

For every $s\in S(\widehat D)$, because $G_s$ is bounded, its
underlying set is contained in $R^n$, for some definable real closed
field $R$ in $\widehat D$. The field $R$, with all its $\widehat
D$-induced structure is o-minimal, and therefore, Corollary
\ref{cor1.2} applies and hence $Z(G_s)\cap [G_s,G_s]_n$ is finite. Since
this is true for every $s\in S(\widehat D)$, there exists a bound
$k=k(n)$ such that $\widehat D\models \forall s\in S \,\,|Z(G_s)\cap
[G_s,G_s]_n|\leq k.$ This is a first order statement which carries
over to $\widehat \CN$ and therefore to $\CN$ as well. It follows
that $Z(G)\cap [G,G]_n$ is finite.\qed
\\

\noindent
{\bf Question.} Is there a direct proof, avoiding the  Lie
algebra argument for the following result: If $\tilde G$ is a
definably connected central extension of a semisimple group in an
arbitrary o-minimal structure then the set $Z(\tilde G)\cap [\tilde
G,\tilde G]_n$ is finite for every $n$?
\\
%{\bf Udi asked here something about central extensions constructed
%by Steinberg but it was not clear to me (or Anand) what this was
%referring to}

\section{The main results}

\subsection{Interpretability results}

\begin{theorem}\vlabel{thm1} Let $\CM$ be an o-minimal expansion of
a real closed field and let $A\to \tilde G\to G$ be a definable
central extension of a  semisimple group $G$, with $\tilde G$
definably connected. Then
\begin{enumerate}
\item $\la A, \tilde G, i,\pi \ra$ is bi-interpretable, over an
imaginary  parameter $\bar c $, with $\la G, A \ra$. The
parameter $\bar c$ names a map from a family of definably connected
components of a $0$-definable set in $\la G,\cdot\ra$ onto
a finite subset of $A$.
\item The exact sequence $A\to \tilde G\to G$ is elementarily
equivalent,  after naming parameters on both sides, to a
semialgebraic central extension $A'\to \tilde G'\to G'$, defined
over the real algebraic numbers, with $\dim(\tilde G')=\dim (\tilde G)$.
\end{enumerate}
 If
$\tilde G$ is definably compact then, in both (1) and (2),
 it is sufficient to assume that $\CM$ expands an ordered group.
In this case $\tilde G'$ of (2) can be chosen definably compact as
well.

\end{theorem}

An immediate corollary of (2) above is:
\begin{cor}\vlabel{compact1} If $\CM$ is an o-minimal expansion of an
ordered group then every definably compact, definably connected
group  is elementarily equivalent, in pure group language, to a
compact semialgebraic (in particular Real Lie) group over the real
algebraic numbers.
\end{cor}

\noindent {\em Proof of Theorem \ref{thm1}}:

(1) We need only to establish the three assumptions of Theorem
\ref{interpret}.

By Theorem \ref{T1.3}, $\la G,\cdot\ra$ has property $\rho$ with
respect to the pure group structure. By Corollary \ref{cor1.2} (and
Corollary \ref{cor1.4} in the definably compact case), $\tilde G$
has the property $Z$. By Claim \ref{perfect1}, $G$ is perfect.\qed

(2) By (1), $\la A,\tilde G, G,i,\pi\ra$  is bi-interpretable with
$\la G, A\ra$, after naming the necessary map. By \ref{T1.3}, $G$
itself is definably isomorphic in $\CM$ to a semialgebraic group
$G'$ over the real algebraic numbers, which is clearly definably
compact if $G$ is. Therefore, in order to prove (2), it is
sufficient to show: Given a finite set $C\sub A$, the structure $\la
A,+,C\ra$ is elementarily equivalent to $\la A',+,C'\ra$ for some
semialgebraic group $A'$ over the real algebraic numbers, and a
finite subset $C'\sub A'$.  This is exactly the content of Claim
\ref{abelianclaim} in the Appendix.  \qed

\noindent{\bf Remark} As the proof of (2) above shows,  the only
obstacle for a definable central definably connected extension
$\tilde G$ of a definable  semisimple group to be definably\emph{
isomorphic} to a semialgebraic group is the group $Z(\tilde G)$.
Hence, if $Z(\tilde G)$ is definably isomorphic in $\CM$ to a
semialgebraic group then so is $G$. We can also prove an analogue
for algebraic groups:
\\

\begin{cor}\vlabel{alg} Assume that $\CM$ is an o-minimal expansion of a real closed
field $R$.
 Let $A\to \tilde G\to G$ be a definable central extension of a definably connected
semisimple group $G$.

If $G$ is a stable group and $A$ is definably isomorphic in $\CM$ to
an algebraic group over $K=R(\sqrt{-1})$, then $\tilde G$ is
definably isomorphic in $\CM$ to an algebraic group over $K$.
\end{cor}
\proof Since $G$ is stable, $G/Z(G)$ is a direct product of
definably simple stable groups, which we may assume are all
algebraic groups over $K$ (see \ref{facts}). Hence, $G/Z(G)$ is
definably isomorphic to an algebraic group over $K$. By Theorem
\ref{T1.3}(1), $G$ is definable, possibly over parameters, in
the group  $G/Z(G)$ and therefore it is definably isomorphic in
$\CM$ to algebraic group over $K$. We continue as in the proof
of Theorem \ref{thm1}.  \qed
\\

We end this discussion with an example showing that not every definably connected
group in an o-minimal structure is elementarily equivalent to a real Lie group which is
definable in an o-minimal structure.
This is a small variation of an example from \cite{pps3}, so we will be brief:
\\

\noindent{\bf Example} Let $\CM=\la R,<,+,\cdot,exp\ra$ be a nonstandard model of theory of
the real exponential field, and let $\alpha\in R$ be element greater than all natural numbers.
We define
$$G=\left\{ \left(
                            \begin{array}{ccc}
                              t & 0 & u \\
                              0 & t^{\alpha} & v \\
                              0 & 0 & 1\\
                            \end{array}
                          \right): u,v\in M, t>0\right\}.$$

The group $G$ is a solvable centerless group, and as is shown in \cite{pps3}*{p.4}, the structure
$\CM_{\alpha}=\la R,<,+,\cdot,t\mapsto t^{\alpha}\ra$ is interpretable in the pure group $G$. If $G$ were
elementarily equivalent to a definable
real Lie group $H$ in some o-minimal structure over the reals then $H$ would interpret a structure
$\CN_{\alpha}\equiv \CM_{\alpha}$ so the underlying field of $\CN_{\alpha}$ is non-archimedean. However,
every real closed field which is interpretable in an o-minimal structure over the reals must be archimedean
(its ordering is Dedekind complete). Contradiction.
\\

\subsection{Structural results}
 We can now deduce a structural result about
definably compact groups in o-minimal structures.

\begin{cor}\vlabel{structure} Let $\CM$ be an o-minimal expansion of an ordered
group. Let $\tilde G$  be a definably compact, definably connected
group.
 Then
 \\(i) $[\tilde G,\tilde G]$ is definable, definably connected and semisimple.
 \\(ii)  $\tilde G$
equals the almost direct product of $[\tilde G,\tilde G]$ and
$Z(\tilde G)^0$. Namely, $\tilde G$ is the group product of $[\tilde
G,\tilde G]$ and $Z(\tilde G)^0$ and the intersection of two groups
is finite. In particular, $\tilde G\simeq(Z(\tilde G)^0\times
[\tilde G,\tilde G])/F$, for a finite central subgroup $F\sub
Z(\tilde G)^0\times [\tilde G,\tilde G]$.
\end{cor}
\proof (i) By \ref{compact1}, $\tilde G$ is elementarily equivalent
to a compact real Lie group $H$. By classical Lie group theory,
$[H,H]$ is a closed connected semisimple subgroup of $H$ (Indeed,
this can be found, for example, in \cite{OV}*{Chapter 5.2, Theorem
4}).

 By topological
compactness, there exists a $k$ such that $[H,H]_k=H$, i.e. the set
$[H,H]_k$ is already a subgroup of $H$.
 It follows that the same is true for $\tilde G$ hence $[\tilde G,\tilde G]=[\tilde G,\tilde G]_k$
  is definable.
The group $[\tilde G,\tilde G]$ is definably connected as the
continuous image of the definably connected group $\tilde G$.

(ii) Because the intersection of $Z(\tilde G)$ and $[\tilde G,\tilde
G]$ is 0-dimensional (Corollary \ref{cor1.2}), it must be finite.
the final clause is immediate from the fact that $[\tilde G,\tilde
G]$ acts trivially on $Z(\tilde G)^0$, and hence the map
$(g,h)\mapsto gh$ from $Z(\tilde G)^0\times [\tilde G,\tilde G]$ to
$\tilde G$ is a homomorphism with finite kernel.\qed

\noindent{\bf Remark } Actually, by a theorem of Goto,
\cite{HM}*{Theorem 6.55}, every element of $[H,H]$ is a commutator
(i,e, $[H,H]=[H,H]_1$) hence the same is true in every definably
compact group.

\begin{cor}\vlabel{cor1.5}  Let $\CM$ be an o-minimal expansion of an ordered
group and assume that $\tilde G$ is definably connected with $\tilde
G/Z(\tilde G)$  semisimple and definably compact. Then \\(i)
$[\tilde G,\tilde G]$
 is definable.
 \\(ii)  $\tilde G\simeq (Z(\tilde G)^0\times
 [\tilde G,\tilde G])/F$, for a finite central subgroup $F\sub Z(\tilde G)^0\times
 [\tilde G,\tilde G]$
 \end{cor}
 \proof By \ref{thm1} (2), the assumption implies that $Z(G)^0\to \tilde G\to G$
 is elementarily equivalent to a semialgebraic extension $Z_0\to
 \tilde G_0\to G_0$ over the real numbers, with $G_0$ a compact connected semisimple
 real  Lie
 group. By Corollary \ref{cor1.2} (ii), we have $L(\tilde G_0)\simeq L(Z_0)\oplus [L (\tilde
 G_0), L (\tilde G_0)]$, and the Lie algebra of
 the locally definable group $[\tilde G_0,\tilde G_0]$ equals
 $[L(\tilde G_0),L(\tilde G_0)]$. Moreover, $[L(\tilde G_0),L(\tilde G_0)]$
 is isomorphic to $L(G_0)$.

 We first claim that $[\tilde G_0,\tilde G_0]$ is a compact
subgroup of $\tilde G_0$.

We  recall the following definition: A Lie algebra over $\reals$ is
called {\em compact} if it admits an invariant positive definite
scalar product.
 Clearly, a subalgebra of a compact Lie algebra is also compact and
 if a Lie algebra is commutative then it is compact (any positive definite scalar product
 will do). Furthermore, the direct sum
 of two compact Lie algebras is compact as well.
 Finally, the Lie algebra of a compact real Lie group is a compact, see \cite{OV}*{p.  228, 12}.

Because $G_0$ is compact, its Lie algebra is compact. Hence
$L(G_0)=[L(\tilde G_0),L(\tilde G_0)]$ is a compact Lie algebra as
well.
 Because $L(Z_0)$ is abelian it is also compact. It follows
 that $L(\tilde G_0)$ is a compact Lie algebra as well.
 We now apply a theorem about connected Lie groups  with compact Lie algebras
  (see \cite{OV}*{p.  242, Theorem 5})
 and conclude that, as a Lie group,
 $\tilde G_0=B\times C[\tilde G_0, \tilde G_0]$, for Lie subgroups
 $B,C\sub Z(\tilde G)^0$ ($B$ torsion-free and $C$ compact), and with
 $C\cdot [\tilde G_0,\tilde G_0]$ a compact Lie subgroup, which we
 denote  by $H$ (these groups are not claimed to be definable).

Because $L([\tilde G_0,\tilde G_0])$ is a maximal semisimple
subalgebra of $L(H)$ it follows, using Levi decomposition theorem,
that $L(H)=L(C)\oplus L([\tilde G_0,\tilde G_0])$. As we saw before,
since $H$ is compact, the group $[H,H]$ is a compact subgroup of $H$, which in this
case must equal $[\tilde G_0,\tilde G_0]$. Thus, $[\tilde G_0,\tilde
G_0]$ is a compact subgroup of $\tilde G_0$.

 It follows that there is a $k$ (actually, as remarked above, $k=1$),
such that $[\tilde G_0,\tilde G_0]_k=[\tilde G_0,\tilde G_0]$. By
elementary equivalence, $[\tilde G,\tilde G]_k=[\tilde G,\tilde G]$,
hence this group is definable.

As we already saw, see \ref{cor1.2}, the intersection of $[\tilde
G,\tilde G]$ and $Z(\tilde G)$ has zero dimension and therefore is
finite. Hence, $\tilde G\simeq (Z(\tilde G^0)\times [\tilde G,\tilde
G])/F$ for a finite central subgroup $F$.\qed

\section{The connection to Pillay's conjecture}

Corollary \ref{compact1}, gives a strong connection between
definably compact, definably connected groups in o-minimal
structures and compact real Lie groups.

Pillay's Conjecture (now a theorem, for expansions of ordered
groups, see \cite{nip1} and \cite{nip2}) suggests another such
connection to compact real Lie groups:
\\

{\em Let $G$ be a definably compact group in a sufficiently
saturated o-minimal structure. There exists a minimal type-definable
subgroup $G^{00}\sub G$ such that $G/G^{00}$, equipped with the
logic topology, is isomorphic to a real Lie group, and the
topological dimension of $G/G^{00}$ equals the o-minimal dimension
of $G$.}
\\

%{\bf Udi asked here two questions: 1. Is it enough to assume a local
%group for PC? 2. Is it true that every definable group lives in the
%finite union of group-intervals? I don't know the answer to either
%of those. I think that a positive answer to (2) could be very nice.
%I cannot see an immediate way to get it. It is probably a good
%project for someone, maybe a student}

Our goal is to prove, in the definably connected case, that the pure
groups $G$ and $G/G^{00}$ are elementarily equivalent. More
precisely, we will prove:
\begin{theorem}\vlabel{thm4} Let $\CM$ be a sufficiently saturated
 o-minimal expansion of
an ordered group and let $G$ be a definably compact, definably
connected group. Then $$\la G, \cdot \ra \equiv \la
G/G^{00},\cdot\ra.$$

Moreover, the map $\pi: G\to G/G^{00}$ ``splits elementarily'',
namely there exists an elementary embedding (with respect to the
group structure)  $\sigma:G/G^{00}\to G$ which is also a section for
$\pi$.
\end{theorem}

By Corollary \ref{structure}, definably compact groups can be
analyzed using abelian and semisimple subgroups. We first handle
the abelian case. Because we are going later on to treat definable
groups which are not necessarily connected, we will need to work in
a more general setting of abelian groups together with finitely many
automorphisms.

 \subsection{Definable abelian
groups with an additional abelian structure}

Let $A$ be an abelian group definable in an o-minimal structure
$\CM$, which we assume to expand an ordered group. We denote by
${\bf A}_{ab}$ the expansion of the group $A$ by all $\CM$-definable
subgroups of $A^n$, $n\in \mathbb N$, and let $L_{ab}$ be the
associated language (note that if $A$ is definably compact and
abelian then, by \cite{ps1}*{Cor. 5.2}, every
 $\CM$-definable subgroup of $A$ is actually $0$-definable in $\CM$).
For $B$ a subgroup of $A$ we let ${\bf B}_{ab}$ be the
$L_{ab}$-substructure of ${\bf A}_{ab}$.

 In the appendix we treat the general (not necessarily o-minimal)  such
situation and observe, using known results:

\begin{fact}\vlabel{module1}
Let $A$ be a abelian definable group in an o-minimal
structure, then
\begin{enumerate}
\item The structure ${\bf A}_{ab}$ eliminates quantifiers.

\item Assume that $B\leq A$ is an arbitrary subgroup of $A$.

Then ${\bf B}_{ab}\prec {\bf A}_{ab}$ if and only if the following
hold: \\(i) For every $0$-definable (in ${\bf A}_{ab}$)
subgroup $S\leq A^{n+k}$ and $b\in B^k$,
$$S(B^n,b)\neq \emptyset \Leftrightarrow S(A^n,b)\neq
\emptyset$$ and (ii) For every $0$-definable (in ${\bf
A}_{ab}$) subgroups $S_1 \leq S_2 \leq A^n$,
$$[S_2:S_1]= [S_2\cap B^n:S_1\cap B^n],$$ with the meaning that if
this index is infinite on one side  then it is infinite on the
other.
\item If ${\bf B}_{ab}\prec{\bf A}_{ab}$, for a subgroup $B$ of $A$, then there exists a surjective
group homomorphism $\phi:A\to B$ which is the identity map on $B$
and in addition sends every $0$-definable group $S\sub A^n$ onto
$S\cap B^n$. (We call such a map $\phi$ {\em a homomorphic
retract}).
\end{enumerate}
\end{fact}

\begin{lemma}\vlabel{Lab}
Let $A$ be a definable compact, definably connected abelian group,
in an o-minimal expansion $\CM$ of an ordered group and let ${\bf
A}_{ab}$ be as above.
\begin{enumerate}

\item  Assume that
${\bf B}_{ab}\prec {\bf A}_{ab}$ for a subgroup $B$ of $A$,
$\phi:A\to B$ a homomorphic retract and let $A_1=ker \phi$. Then
$({\bf Tor(A)+A_1})_{ab}\prec {\bf A}_{ab}$.

\item If ${\bf A}_{ab}$ is sufficiently saturated then $({\bf Tor(A)+A^{00}})_{ab}\prec {\bf
A}_{ab}$.
\end{enumerate}
\end{lemma}

Note that if we take, in (1) above, $B=A$ and $\phi$ the identity
map then the lemma implies in particular that ${\bf
Tor(A)}_{ab}\prec {\bf A}_{ab}$.

\proof (1)
 We need to see that $Tor(A)+A_1$ satisfies the two
requirements of \ref{module1}(2). Note that $B$ is divisible and
contains all torsion elements of $A$ (there are finitely many for
each exponent). Therefore, since  $A=B\oplus A_1$, the group $A_1$
is divisible as well (see \ref{abelian1}).

Clause (i): Let $S$ be an $\CM$-definable subgroup of $A^{n+k}$ and
assume that $(b+a_1,c)\in S$ for some $b\in B^n, a_1\in A_1^n$ and
$c\in (Tor(A)+A_1)^k$. We want to show that there is $a\in
(Tor(A)+A_1)^n$ such that $(a,c)\in S$.

We write $c=c_1+c_2$, for $c_1\in Tor(A)^k$ and $c_2\in A_1^k$.
 If $mc_1=0$ then
 $$m(b+a_1,c)=m(b+a_1,c_1+c_2)=(mb+ma_1,mc_2)\in S,$$  with $(ma_1,mc_2)\in
 A_1^{n+k}$. Because $\phi$ is a retract, we have
 $(mb,0)\in S$. Now, if $k$ is the index of $S(A^n,0)^0$ in
 $S(A^n,0)$ then $(kmb,0)\in S(A^n,0)^0$. Because this last group is divisible,
 there exists $(a,0)\in S(A^n,0)^0$ such
 that $kma=kmb$ and therefore $b-a\in Tor(A)^n$. Finally, we have
 $(b+a_1,c)-(a,0)=(b-a+a_1,c)\in S$, with $b-a+a_1\in (Tor(A)+A_1)^n$,
 as needed.

Clause (ii): We will actually prove a stronger statement than needed
here:

 {\em (*) If $C\sub A^n$ is a divisible subgroup containing $Tor(A)^n$
then for every $S_1\sub S_2\sub A^n$ definable groups.
$$[S_2:S_1]=[S_2\cap C:S_1\cap C]$$ (and both are
infinite if one of them is)}.\\

We clearly have $[S_2\cap C:S_1\cap C]\leq [S_2:S_1]$, so we only
need the opposite inequality.

 Assume first that $S_1$ is definably connected, hence divisible.
 It follows that every torsion
element of the group $S_2/S_1$ (i.e. a coset of $S_1$) contains a
torsion element of $S_2$ and therefore an element of $S_2\cap C$.
Hence, we have an injective map from $Tor(S_2/S_1)$ into $S_2\cap
C/S_1\cap C$ and hence $|Tor(S_2/S_1)|\leq [S_2\cap C:S_1\cap C]$.

 If $[S_2:S_1]$ is infinite then $dim S_1<dim S_2$ and
therefore $S_2/S_1$ is a definably compact group of positive
dimension. It follows from \cite{edot} that $Tor(S_2/S_1)$ is
infinite and therefore, by the inequality above, so is $[S_2\cap
C:S_1\cap C]$. If $S_2/S_1$ is finite then all its elements are
torsion and therefore, by the same inequality we have
$$S_2/S_1=Tor(S_2/S_1)=[S_2\cap C:S_1\cap C]$$

 If $S_1$ is not definably connected
then we apply the above argument first to $[S_2:S_1^0]$ and
$[S_1:S_1^0]$ and then conclude the result for $[S_2:S_1]$.

(2) Since $A^{00}$ is divisible we only need, by (*) above,
 to see that Clause (i) holds for $Tor(A)+A^{00}$.

Let $S\sub A^{n+k}$ be a $\CM$-definable group, $\pi_2:S\to A^k$ the
projection map onto the second coordinate and let $S_1=\pi_2(S)$. As
we saw in \cite{nip1}, we have
$$S^{00}=(A^{n+k})^{00}\cap S=(A^{00})^{k+n}\cap S$$ and
$$\pi_2(S^{00})=S_1^{00}=(A^{00})^k\cap S_1.$$

Assume now that $(a,c_1+c_2)\in S$ for $a\in A$, $c_1\in Tor(A)$,
$c_2\in A^{00}$. If $mc_1=0$ then $(ma,mc_2)\in S$, with $mc_2\in
S_1^{00}$. Because $S_1^{00}=\pi_2(S^{00})$ there exists $e_1\in
(A^{00})^n$ such that $(e_1,mc_2)\in S$. Moreover, because $A^{00}$
is divisible, we have $e_1=me$ for some $e\in (A^{00})^n.$ It
follows that $(ma,mc_2)-(me,mc_2)=(ma-me,0)\in S$.

If we let $k$ be the index of  $S(A^n,0)^0$ in $S(A^n,0)$ then
$(kma-kme,0)\in S(A^n,0)^0$ and there exists $(d,0)\in S$ such that
$kmd=kma-kme$. In particular, $(a-e)-d\in Tor(A)^n$. We now have
$(a,c_1+c_2)-(d,0)=(a-d,c_1+c_2) \in S$, with $a-d\in
(Tor(A)^n+A^{00})^n,$ as needed.\qed

By considering the special case of a compact real Lie group
definable in the o-minimal structure $\reals_{an}$ (or by a direct
modified version of the above proof) we also have:
\begin{lemma}\vlabel{compact-ab} Let $B$ be a connected, compact abelian real Lie
group and let ${\bf B}_{an}$ be the expansion of $(B,+)$ by adding a
predicate for every compact Lie subgroup of $B^n$. Then
$${\bf Tor(B)}_{an}\prec {\bf B}_{an}.$$
\end{lemma}

We can now state the main result in the abelian case:

\begin{theorem}\vlabel{theorem-ab}
Let $A$ be a definably compact, definably connected abelian group in
a sufficiently saturated  o-minimal expansion $\CM$ of an ordered
group. Endow $B=A/A^{00}$ with an $L_{ab}$-structure by interpreting
$R_S$, for every $0$-definable subgroup $S\sub A^n$, as the group
$\pi(S)\sub B^n$ (where $\pi:A\to A/A^{00}$ is the projection map).
Let ${\bf B}_{ab}$ be the induced structure on $B$. Then,
\begin{enumerate} \item  ${\bf B}_{ab}\equiv {\bf A}_{ab}$.
Moreover, there is an $L_{ab}$-elementary embedding $\sigma:B\to A$
which is a section for $\pi$.
\item The structure ${\bf B}_{ab}$ is a reduct of ${\bf B}_{an}$ above.
\end{enumerate}
\end{theorem}
\proof (1) We start with the structure $T=({\bf
Tor(A)+A^{00}})_{ab}$, which, by \ref{Lab} (2), is an elementary
substructure of ${\bf A}_{ab}$. By \ref{module1}(3), there exists a
homomorphic retract $\phi:A\to T$ and if we let $A_1=ker \phi$ then,
by \ref{Lab}(1), the structure ${\bf C}_{ab}=(Tor(A)+A_1)_{ab}$ is
also elementary in ${\bf A}_{ab}$.

It is left to see that  the restriction of $\pi$ to ${\bf C}_{ab}$
induces an isomorphism of ${\bf C}_{ab}$ and ${\bf B}_{ab}$.

Let $S\sub A^n$ be an $\CM$-definable group and $c\in C^n$. We need
to see that $\pi(c)\in \pi(S)$ if and only if $c\in S$. Write
$c=a+a_1$ for $a\in Tor(A)^n, a_1\in A_1^n$, and assume that
$\pi(a+a_1)\in \pi(S)$. It follows ($A^{00}\in ker\pi$) that for
some $b\in (A^{00})^n$ we have $a+b+a_1\in S$. Because $A_1$ is the
kernel of the retract $\phi$, we have $a+b\in S$. Because $a$ is a
torsion element, there exists an $m$ such that $$ma+mb=mb\in S\cap
(A^{00})^n=S^{00}.$$ Because $S^{00}$ is divisible, there exists
$b_1\in S^{00}$ such that $mb_1=mb$ and therefore $b-b_1$ is a
torsion element. However, $b-b_1$ belongs to the torsion-free group
$(A^{00})^n$, hence $b=b_1$ and  $b\in S$. We can therefore conclude
that $c=a+a_1\in S$, and therefore $\phi|C$ is an isomorphism of
${\bf C}_{ab}$ and ${\bf B}_{ab}$. The inverse map $\sigma:B\to C$
is an elementary embedding of ${\bf B}_{ab}$ into ${\bf A}_{ab}$.

(2) The image under the projection map of every definable set in
$A^n$ is closed in the Euclidean topology on $(A/A^{00})^n$. Since
every closed subgroup of  Lie group is itself a Lie subgroup, it
follows that for every definable $S\leq A^n$, $\pi(S)$ is a Lie
subgroup of $B^n$. Therefore, ${\bf B}_{ab}$ is a reduct of ${\bf
B}_{an}$.
 \qed
\\

%\begin{lemma} \vlabel{automorab} Let $G$ be a definably compact
%(not necessarily connected) group in a sufficiently saturated
%o-minimal structure, with $G^0$ abelian, and let $\pi:G\to G/G^{00}$
%be the projection map. Then $G$ and $G/G^{00}$ are elementarily
%equivalent in the pure group language. Moreover, there exists an
%elementary embedding $\sigma:G/G^{00}\to G$ which is a section for
%$\pi$.
%\end{lemma}
%\proof Let $A=G^0$, and $\widehat A=A/A^{00}$.

%Then $G=(A\rtimes F)/K$ for a finite group $F$ acting definably on
%$A$ and a finite normal subgroup K. As in \ref{ }, consider the
%structure ${\bf A}_{ab}$ which $\CM$ induces on $A$ by adding a
%predicate to each $\CM$-definable subgroup of $A^n$, $n\in \mathbb
%N$ and let $\widehat {\bf A}_{ab}$ be the expansion of $\widehat A$
%by the image of every definable subgroup of $A^n$. By \ref{ }, there
%is an elementary embedding of $\widehat {\bf A}_{ab}$

\noindent
{\bf Remark.} Note that for all the results above we did
not require the full theorem of Edmundo-Otero about the structure of
$Tor(A)$ for a definably compact abelian group $A$. We only needed
the weaker statement that every definably compact infinite group has
infinitely many torsion elements. However, without the stronger
result we will not be able to conclude that $\dim(A/A^{00})=\dim A$.
\\

\noindent
{\bf Remark.} It is not hard to see that ${\bf B}_{ab}$ is $\omega$-saturated, hence it and ${\bf A}_{ab}$
will actually be
$L_{\infty, \omega}$-equivalent, improving 7.5(1).
\\

\subsection{The general case}

\noindent{\bf Proof of Theorem \ref{thm4}}: By \ref{structure} and
our earlier analysis, the group $H=[G,G]$ is definable, semisimple,
$Z(G)\cap H$ is finite. We have $G\simeq (Z(G)\times H)/F$, for a
finite central subgroup $F\sub Z(G)\times H$. Because $G$ is
definably connected, $H$ is definably connected as well. Hence, we
also have $G\simeq (Z(G)^0\times H)/F_1$, for a finite central
subgroup $F_1$. We write $A=Z(G)^0$.

For any definably compact, definably connected  $K$, we denote by
$\widehat K$ the group $K/K^{00}$. For $a\in K$, we let $\hat
a=\pi(a)\in \widehat K$.

We prove Theorem \ref{thm4} in several steps.

\noindent {\bf Claim I } If $A$ is abelian and definably connected
then for every finite subgroup $A_1\sub A$,
$$\la A,\cdot,\{a:a\in A_1\}\ra\equiv \la \widehat {A},\cdot,\{\hat
a:a\in A_1\}\ra.$$
 Moreover, there exists an elementary section
$\sigma_A:\widehat A\to A$ with $\sigma_A(\hat a)=a$ for every $a\in
A_1$.
\\

\proof Because $A_1$ is $\CM$-definable, this is almost immediate
from \ref{theorem-ab}. We only need to notice that since all
elements of $A_1$ are torsion elements, the elementary embedding of
$\widehat A$ into $A$ necessarily sends every $\hat a\in \widehat
A$, $a\in A_1$, to the element $a$.

\noindent{\bf Claim II} Given $H$ definably connected,
 definably compact and semisimple, for every finite central subgroup $H_1\sub
H$,
$$\la H,\cdot,\{a:a\in H_1\}\ra\equiv \la \widehat {H},\cdot,\{\hat
a:a\in H_1\}\ra.$$ Moreover, there exists an elementary
embedding $\sigma_H:\widehat H\to H$ which is a section for $\pi$,
such that $\sigma_H(\hat h)=h$ for every $h\in H_1$.
\\

\proof  By \ref{T1.3}, we may assume  that $H$ is a semialgebraic
group definable over the real algebraic numbers. Hence, every
$0$-definable set contains an element from
 $dcl(\emptyset)$, and moreover every finite group is $0$-definable.
  We thus may assume
that $H$ is semialgebraic and definable over the real algebraic
numbers.

In this case (see \cite{nip1} for details), there exists an elementary embedding
$\sigma_H: \hat H\to H$ which is also a section for $\pi:H\to \hat
H$, and in particular, $H=H^{00} \rtimes \sigma_H(\hat H)$. Because
$\sigma_H$ is elementary, for every $\hat h\in \widehat H_1$, the
element $\sigma_H(\hat h)$ is a torsion element of $H$ and we have
$\sigma_H(\hat h)\hat h^{-1}\in H^{00}$. However $H^{00}$ is
torsion-free (see \cite{b2}*{Theorem 4.6}) and therefore $\sigma_H(\hat h)=h$.
 \\

\noindent{\bf Claim III} If $G=(A\times H)/F$, for $A$ definably
connected, definably compact abelian, $H$ definably connected,
definably compact semisimple and $F$ a finite central subgroup of
$A\times H$, then $G\equiv \widehat G$. Moreover, there exists an
elementary embedding $\sigma_G:\widehat G\to G$ which is a section
for $\pi$.
\\

\proof Since $F$ is finite it is contained in $A_1\times H_1$ for some
finite groups $A_1\sub A$ and $H_1\sub H$. It is easy to see that
$\widehat{A\times H}\simeq \widehat {A}\times \widehat {H}$. By step
I and step II, we have an elementary section:
$$\sigma:\la \widehat{
{A}}\sqcup \widehat {H},\{\hat a\in \widehat{A_1}\},\{\hat h\in
\widehat {H_1}\}\ra \, \to \, \la A\sqcup H, \{a\in A_1\},\{h\in
H_1\}\ra$$ sending each $\hat a$ and $\hat h$ to $a$ and $h$ respectively. It follows that we have
elementary section
$$ \sigma_1: \la \widehat
{A}\times \widehat {H}, \cdot, \{(\hat a,\hat h)\in \widehat
{A_1}\times \widehat {H_1}\}\ra\, \to \, \la A\times H,\cdot,
\{(a,h)\in A_1\times H_1\}\ra, $$ and hence also
$$\sigma_2 :\la \widehat
{A}\times \widehat {H},\cdot, \{\hat g\in \widehat {F_1}\}\ra\simeq
\la \widehat{A\times H}, \cdot, \{\hat g\in \widehat
{F_1}\}\ra\,\to\, \la A\times H,\cdot, \{g\in F_1\}\ra.$$
This last section sends each $\hat g\in \hat F$ to $g\in F$.

In order to complete the proof of Claim III, it is therefore
sufficient to prove the following general fact (with $K$ now playing
the role of $A\times H$):
\\

\noindent{\bf Fact } Let $K$ be a definably connected, definably
compact group and $F\sub K$ a finite central subgroup. Assume that
$\sigma_K:\widehat K\to K$ is an elementary section of $\pi_K:K\to
\widehat K$ which, for every $g\in F$ sends $\hat g\in \widehat F$
to $g\in F$. Then the map $\sigma_{K/F}$ which sends the element $(gF)(K/F)^{00}$
of $\widehat{G/F}$ to $(\sigma_K(gK^{00}))F\in K/F$ is an elementary section
for $\pi_{K/F}:K/F\to \widehat{K/F}.$
\\

 \proof It is easy to see that the map $\sigma:\widehat K/\widehat F\to K/F$ which sends
 $(gK^{00})\widehat F$ to $\sigma_K(gK^{00})F$ is elementary. It is also not hard to see that the map
 $\sigma': (gF)(K/F)^{00})\mapsto (gK^{00})\widehat F$ is an isomorphism of $\widehat{K/F}$ and $\widehat
 K/\widehat F$ (we use here the fact that
 the projection map $\pi_F:K\to K/F$ sends $K^{00}$ onto $(K/F)^{00}$).

 The composition of $\sigma $ and $\sigma'$ gives the desired $\sigma_{K/F}$.

 %Consider the
%surjective homomorphism.
%$$\pi:K\twoheadrightarrow \widehat{K} \mbox{ with } ker (\pi)=K^{00}.$$
%Because $\pi(F)=\widehat{F}$ we have $K/(K^{00}F)\simeq
%\widehat{K}/\widehat {F}$.
% Similarly, we have
%$$\alpha:K\twoheadrightarrow K/F \mbox{ with } ker(\alpha)=F.$$
%Because $\alpha(K^{00})=(K/F)^{00}$ (it is always true for a
%surjective $\phi:G\to H$ that $\phi(G^{00})=H^{00}$) we have
%$K/(K^{00}F)\simeq \widehat{K/F}$.

%Summing up, we have $\widehat{K/F}\simeq
%\widehat{K}/\widehat{F}$.\qed
%\\

\section{Finite extensions of o-minimal groups}

In this section we consider finite (but not necessarily central)
extensions of arbitrary groups definable in o-minimal structure.
Finite extensions of groups in o-minimal structures were studied by Edmundo, Jones and Peatfield in
\cite{ejp}. The following was shown there,
using universal covers: If $G$ is a definable,  definably
connected group in an o-minimal structure $\CM$ expanding a real closed field,
  and  if $\pi:\tilde G\to G$ is any finite
 definable extension of $G$, defined possibly in an o-minimal expansion $\CN$ of $\CM$,
 then $\tilde G$ is definably isomorphic in $\CN$ to a group definable in $\CM$.

 Indeed, the above result is not stated as such but can be read off the
 proof of Proposition 3.2 in \cite{ejp}.
  This implies for example that if $G$ is semialgebraic then so is
  every finite extension of $G$ which is definable in $\CN$.

In this section we give two different proofs for similar results about
the interpretability  of finite extensions of definable groups and arbitrary topological covers of definable groups over the reals. Although the two main results, Theorem \ref{interpret2} and Theorem \ref{real-interpret2} overlap in the case of finite covers, the assumptions and techniques are different so we include both.

 We first need the following fact about the
structure of arbitrary definably connected groups in o-minimal
structures. For $G$ a group and $n\in \mathbb N$ we let
$\sigma_n:G\to G$ be the map $\sigma_n(g)=g^n$.

 \begin{lemma}\vlabel{div1} Let $G$ be a definably connected group in an o-minimal
 structure. Then,
 \\(i)  The group $G/[G,G]$ is divisible, namely, for every $n\geq 1$,
 $\sigma_n(G)[G,G]=G$. In fact,
  there exists $k\in \mathbb N$ such that $\sigma_n(G)[G,G]_k=G$.
  \\(ii) For every $n$, let $\la \sigma_n(G)\ra $ be the group generated by
  all elements $g^n$, $g\in G$. Then $G=\la \sigma_n(G)\ra$. In fact,
  there is $k\in \mathbb N$ such that $G=\sigma_n(G)\cdots \sigma_n(G)$ ($k$-times).
 \end{lemma}
 \proof (i) We use induction on $\dim G$. If $\dim G=1$ then $G$ is abelian
 and therefore divisible.  We consider the general case.

 Assume first that $G$ has an infinite definable normal abelian subgroup $A$. By induction,
 $G/A$ satisfies the lemma and therefore, for every $n\in \mathbb N$, $A\sigma_n(G)[G,G]=G$.
 Because $A$ is divisible, it is contained in $\sigma_n(G)$ and therefore $\sigma_n(G)[G,G]=G$, as needed.

 If $G$ has no infinite definable normal abelian subgroup then $G$ is semisimple and therefore, by
 Claim \ref{perfect1}, we have $[G,G]=G$.

For the last clause of (i), we may work in a sufficiently saturated structure, where
the existence of such a $k$ is clear. Once proved there, the same $k$ works
for $G$ in any structure.

(ii) As before we may work in a sufficiently saturated structure. For $G$ abelian
the result is clear since it is divisible.

Assume that $G$ has an infinite definable normal abelian subgroup
$A$. By induction on dimension we have $G/A=\la \sigma_n(G/A)\ra$,
which implies that $G=A\la \sigma_n(G)\ra$. However, since $A$ is
divisible it is contained already in $\sigma_n(G)$ and hence $G=\la
\sigma_n(G)\ra$.

 If $G$ has no infinite definable normal abelian subgroup then it is
 semisimple. Let us see why the theorem is indeed true in this case.

 We first
assume that $G$ is definably simple. If $G$ is not definably compact
then it is abstractly simple (see \ref{facts}(5)). The group $\la
\sigma_n(G)\ra$ is clearly invariant under all automorphisms of $G$
hence normal, so $G=\la \sigma_n(G)\ra$. If $G$ is definably
compact, then by \ref{facts} it is elementarily equivalent to a
simple compact real Lie group $H$. By simplicity, $H=\la
\sigma_n(H)\ra=\bigcup_{k=1}^{\infty} \sigma_n(H)\cdots \sigma_n(H)
\mbox{ ($k$-times)}$. It follows from compactness that for some $k$
we have $H=\sigma_n(H)\cdots \sigma_n(H)$ ($k$-times). This implies
that the same is true for $G$.

 If $G$ is semisimple then $Z(G)$ is finite and we have
 $G/Z(G)= H_1\times\cdots\times H_r$, for $H_i$ definably
 simple. By the above, each $H_i$ satisfies $H_i=\la \sigma_n(H_i)\ra$, and hence we have
 $G=Z(G) \la \sigma_n(G)\ra$, so $\la \sigma_n(G)\ra$ has finite
 index in $G$. However, $\la \sigma_n(G)\ra$ is a countable union of
 definable sets and therefore it follows that $G$ is a countable union of such sets.
 Because of saturation, this implies that $\la \sigma_n(G)\ra $ is actually
 definable (in finitely many steps) and by the definable connectedness
 of $G$ we have $G=\la \sigma_n(G)\ra$.\qed

We prove:
\begin{theorem} \vlabel{finite1} Let $\CM$ be an arbitrary structure, sufficiently saturated,
 and let $\CR$ be a definable o-minimal structure in $\CM$.
  Let $G$ be an $\CR$-definable group, $\tilde G$ an
$\CM$-definable  group and let $\pi:\tilde G \to G$ be a $\CM$-definable
surjective homomorphism with finite kernel $N$.

Then $\tilde G$ is internal to $G$ in the reduct containing $\la
\tilde G,\cdot\ra$, $\la G,\cdot\ra$, $\pi$ and a predicate for
$G^0$ the definably connected component of $G$ (denote this reduct
by $\CM'$).

 More precisely,  $\tilde G$ is in the $\CM'$-definable closure of $G$ and a finite
subset of $\tilde G$.
\end{theorem}
\proof Let $\pi:\tilde G \to G$ be the extension map. We may assume
that $G$ is definably connected (in the sense of $\CR$). Indeed,
since $N$ is finite, $\pi^{-1}(G^0)$ has finite index in $\tilde G$,
so if $F\sub G$ is a finite set such that $G=FG^0$ then $\tilde G$
is in the $\CM$-definable closure of $\pi^{-1}(G^0)$ and the finite
set $\pi^{-1}(F)$. It follows that if $\pi^{-1}(G^0)$ is
$\CM'$-internal to $G$ then so is $\tilde G$.

We may also assume that $\tilde G$ has no $\CM'$-definable subgroups
of finite index. Indeed, if $\tilde G_1\sub \tilde G$ is definable
of finite index then $\pi(\tilde G_1)$ has finite index in $G$, and
therefore (we assume $G$ is definably connected)
 $\pi(\tilde G_1)=G$. This in turn implies that $N$ is not contained in $\tilde G_1$, and
therefore $|N\cap \tilde G_1|<|N|$. Using induction on $|N|$ we
could finish.

The assumption that $\tilde G$ has no $\CM'$-definable subgroups of
finite index implies that $N$ is central in $G$. As we will show,
under these assumptions, $\tilde G$ is in the definable closure of
$G$ and $N$.

Let $n=|N|$.

For $k\in \mathbb N$, let $f_k$ be the term $f_k(x_1,\ldots,
x_k)=x_1^n\cdots x_k^n$, and for a group $H$, let $f_{k,H}:H^k\to H$
be the evaluation of the term in $H$.

Let $\pi:\tilde G^k\to G^k$ be the projection map in each of the
coordinates. Similarly to Claim \ref{factor}, we claim that for
$\bar g_1,\bar g_2\in \tilde G^k$, if $\pi(\bar g_1)=\pi(\bar g_2)$
then $f_{k,\tilde G}(\bar g_1)=f_{k,\tilde G}(\bar g_2)$ (we use the
fact that $N$ is central and for every $h\in N$ we have $h^n=e$).

It now follows that there is a $\CM'$-definable surjective map
$h_k:G^k\to \sigma_n(\tilde G)\cdots \sigma_n(\tilde G)$ ($k$-times)
such that $f_{k,\tilde G}$ factors through $\pi$ and $h_k$.

By Theorem \ref{div1}, we may choose $k$ such that
$G=\sigma_n(G)\cdots \sigma_n(G)$ ($k$-times). Said differently, the
map $f_{k,G}:G^k\to G$ is surjective. It easily follows that $\tilde
G=N h_k(G^k)$.\qed

\begin{theorem}\vlabel{interpret2} Let $\CM$ be an o-minimal
structure and assume that  $1\rightarrow N\rightarrow \tilde
G\rightarrow G\rightarrow 1$ is an $\CM$-definable extension with
$N$ finite and $\tilde G$ definably connected.

Let  $\bf G$ be some expansion of $\la G,\cdot\ra$ with property
$\rho$.

Then $1\rightarrow N\rightarrow \tilde
G\rightarrow G\rightarrow 1$  is definably isomorphic in $\CM$
to an a extension $1\rightarrow N_1\rightarrow \tilde
G' \rightarrow G\rightarrow 1$ which is
definable in $\bf G$ over parameters (with $h_{G}:G\to G$ the
identity map).

The parameters name a bijection between a finite
 collection $\mathcal W $ of $\bf G$-definable sets
and a finite subset of $N$. The collection $\mathcal W$ (but not
necessarily each of its sets) is $0$-definable in the pure group
$G$.
\end{theorem}
\proof Note that because $\tilde G$ is definably connected and $N$
is normal and finite then it is necessarily central. The proof of
this theorem is very similar to that of Theorem \ref{interpret}.
Instead of products of $k$ commutators (i.e. the function $F_{k,G}$)
we use products of $k$-many $n$-powers of $G$ (the function
$f_{k,G}$ defined above), with $n=|N|$. Also, instead of $k_n$ we
use here the function $h_k:G^k\to \tilde G$ defined above and
instead of the set $G(k)$ defined there we use the set
$$G_k=\{\bar g\in G^k: f_{k,G}(\bar g)=e\}$$ and its definably
connected components.

Finally, instead of using the fact there that every element of the
perfect $G$ was a product of $k$ commutators, we use Theorem
\ref{div1} which implies that every element of $G$ is a product of
$k$ $n$th-powers. The other details are identical to the proof of
Theorem \ref{interpret}.\qed
\\

\subsection{The real case} Just like in case of Theorem \ref{real-interpret}, if one works
over the field of real numbers then there is no need to assume that $\tG$ is definable
and we obtain the following version of Theorem \ref{interpret2}:

\begin{theorem}\vlabel{real-interpret2} Let $\CM$ be an o-minimal
structure over the real numbers, $G$ an $\CM$-definable group
and assume that  $E:1\rightarrow N\rightarrow \tilde
G\rightarrow G\rightarrow 1$ is a topological extension with
$N$ finite and $\tilde G$ connected.

Let  $\bf G$ be some expansion of $\la G,\cdot\ra$ with property
$\rho$.

Then $E$ is isomorphic as a topological extension, to an extension
$E':1\rightarrow N_1\rightarrow \tilde G\rightarrow G\rightarrow 1$ definable in the structure ${\bf G}$, (with the isomorphism being the identity on $G$).
\end{theorem}
\proof The arguments as to why this version of \ref{real-interpret2} is true are identical to
those explaining the proof of \ref{real-interpret}. In both cases the only o-minimal facts that
are being used apply to $G$ (rather than $\tG$). \qed
\\

We end this diversion into extensions of definable real Lie groups by considering topological covers and related central extensions. This is really an application of work by Edmundo (\cite{e2}) and Edmundo-Eleftheriou (\cite{e3}) on universal covers and local definability in an $o$-minimal setting, as well as work on definable fundamental groups by Berarducci and Otero 
(\cite{bo}). We include the material because we could not find it precisely stated in the literature. In any case thanks to Edmundo for his explanations to us of results implicit in his work, some of which we repeat in the proof below. 

Let us now set up notation for Theorem 8.5 below. $\CM=\la \reals,<,+,\cdot, \cdots\ra$ will be an $o$-minimal expansion of the real field, and $G$ a definably connected group definable in $\CM$ (so $G$ is what we have called a definable real Lie group). $\tG$ will be the topological universal cover of $G$ (also a connected real Lie group) and $\Gamma$ denotes the kernel of $\tG \to G$, namely the fundamental group $\pi_{1}(G)$ of $G$. So $\Gamma$ is a central discrete closed subgroup of $\tG$. If $f:\Gamma \to A$ is a homomorphism from $\Gamma$ into an  abelian group $A$, we form as usual the group ${\tG}_{A} = \tG \times_{\Gamma}A$, and we have a 
central extension $1 \to A \to {\tG}_{A} \to G \to 1$ of $G$ (as abstract groups). We will refer to locally definable groups, for which the reader can consult \cite{e3}, although we give an explanation inside the proof. 

\begin{theorem}\vlabel{universal} (i) $\tG$ and the covering homomorphism can be realized, even topologically, as a locally definable
group and homomorphism in $\CM$.
\newline
(ii) ${\tG}_{A}$ with its group structure, the extension $1 \to A \to {\tG}_{A} \to G \to 1$, together with a section $G\to {\tG}_{A}$, can be interpreted with parameters in the two sorted structure consisting of $\CM$ and 
$\la A, + \ra$.
\end{theorem}
\proof We will be brief. But note first that taking $A = \Gamma$ and $h$ the identity, (ii) says that $\tG$ can be interpreted in the two sorted structure consisting of $\CM$ and $\la \Gamma, + \ra$.

\vspace{2mm}
\noindent
Recall first that for an arbitrary central group extension $E:1\to A\to H\to_{\pi} G\to 1$,
 if $s:G\to H$ is a section for $\pi$, and $h_s(x,y)=s(xy)^{-1}s(x)s(y)$, then $h_s$ (which is called a 2 co-cycle) is a map from $G\times G\to A$ and the group $H$ is isomorphic to the group $H'$ whose underlying set is 
$G\times A$ and whose group operation is given by $(x,a)\cdot (y,b)=(xy,h_s(x,y)+a+b)$ (we can write the second coordinate additively because $A$ is abelian). Moreover $\pi':H'\to G$ is the usual projection, the embedding of $A$ into $H'$ is given by $a\mapsto (1,-a-h_s(1,1))$, and the section $s$ is just $x\mapsto (x,0)$. Hence, in order to recover $E$ we only need to find such a definable co-cycle $h_s$.

\vspace{2mm}
\noindent
We now prove part (i). The statement of (i) is that
there exists in $\CM$
a locally definable group
$\CU=\bigcup_{i\in I}X_i$ (i.e. a bounded directed union of definable sets) with a locally definable group operation (i.e. it is  definable when restricted to each $X_i\times X_j$), and a locally definable surjective homomorphism  
$w:\CU \to G$, and moreover the group $\CU$ with its topology as a locally definable group in $\CM$ is precisely the universal covering of $G$. The reasoning is as follows: Choose $w:\CU \to G$ to be the universal locally definable cover of $G$ as described in \cite{e3}. $\CU$ has of course a topology as a locally definable group, and as the underlying set of $\CM$ is $\reals$, it will be locally Euclidean, connected, and a topological cover of $G$. So $w$ induces an embedding 
$w_{*}$ of the topological fundamental group $\pi_{1}(\CU)$ of $\CU$ into $\pi_{1}(G)$. We will point out that $\pi_{1}(\CU) = 0$, whereby $\CU$ {\em will be} the universal cover $\tG$ of $G$. 
Let $c\in \pi_{1}(\CU)$. So $w_{*}(c) \in \pi_{1}(G)$. By \cite{bo}, $\pi_{1}(G) = \pi_{1}^{def}(G)$ (the definable fundamental group of $G$), whereby $w_{*}(c)$ is represented by a definable path $\gamma$ (beginning and ending at the identity). By Lemma 2.7(1) of \cite{e3}, $\gamma$ lifts to a definable path $\gamma'$ in $\CU$ starting at the identity. On general topological grounds, $\gamma'$ is a loop, and represents $c$. But $\pi_{def}(\CU)$ is trivial, whereby $\gamma'$ is definably homotopic to the identity. Thus $c = 0$.

\vspace{2mm}
\noindent
For (ii) let us first prove the special case that 
$\CU$ is definable in $\la \CM;\Gamma\ra$.
Because $I$ is bounded, there is an $i_0\in I$ such that $X_{i_0}$ projects onto $G$. Because of definable choice, there is an $\CM$-definable $Y\sub X_{i_0}$ and an $\CM$-definable $s:G\to Y$ which is a section for $w$. Moreover, because of local definability, there is $j\in I$ such that $Y\cdot Y\cdot Y^{-1}\sub X_j$, hence the associated 2 co-cycle $h_s:G\times G\to X_j$ is also $\CM$-definable and its image is contained in $X_j$. Note that since $w|X_j:X_j\to G$ is
definable it follows that the image of $h_s$ in $\Gamma$ must be finite (otherwise the kernel of $\pi$ in $X_j$ will be an infinite definable discrete set). Finally, as mentioned above, given the co-cycle $h_s$ we can recover a definable covering $1\to \Gamma'\to \CU'\to G\to 1$ in $\la \CM;\la \Gamma, +\ra\ra$ which is isomorphic to the original one. Because $\Gamma$ is a bounded set (independently of the model $\CM$) the set $\CU'=\Gamma\times G$ can be written as a directed union of definable sets. Since
 $h_s$ is an $\CM$-definable map the group structure on $\CU'$ is locally definable. Finally, the isomorphism $(x,a)\mapsto x\cdot a$ from $\CU'$ to $\CU$ is locally definable as well and therefore a homeomorphism.
 
 \vspace{2mm}
 \noindent
 We now consider the general case of (ii). By what has been done so far we may identify $\tilde G$ with $\CU$. Let
 $h_{s}:G\times G \to \Gamma$ be the $2$-cocycle from the previous paragraph. Define $h':G\times G \to A$ to be
 $f\circ h_{s}$  (where recall $f$ is the given homomorphism of $\Gamma$ into $A$). Then $h'$ is precisely the $2$-cocycle determining ${\tG}_{A}$. As $h_{s}$ was definable in $\CM$ with finite image in $\Gamma$, it follows that the group
 operation on $G\times A$ given in the first paragraph of the proof is definable in the two sorted structure
 consisting of $\CM$ and $\la A, + \ra$. 
\qed
   
\vspace{5mm}
\noindent
Let us remark in closing this section that Theorem 8.5 gives an interesting twist on certain covering structures
considered by Zilber, such as the two sorted structure  $M_{0}$,
say, consisting of $\la {\mathbb C}, + \ra$ in one sort,
$\la {\mathbb C}, +, \cdot \ra$ in the other sort and the complex exponential map $exp$ going from the first sort to the second. The kernel of $exp$ can be identified with the (definable in $M_{0}$) subgroup $\mathbb Z$ of the first sort. It is easy to see that the full structure $M_{0}$ cannot be interpreted in the  reduct consisting of the sorts $\la {\mathbb Z},+\ra$ and $\la {\mathbb C}, +, \cdot \ra$. But Theorem 8.5 says that $M_{0}$ can be so interpreted if we add a predicate for $\mathbb R$ to the second sort. 

We should also mention the Ph.D. thesis of  Misha Gavrilovich on the model theory of the universal covering spaces of complex algebraic varieties, which contains ideas and constructions related to ours above.

%As we saw in Theorem \ref{T1.3}, every definable semisimple group
%in an o-minimal structure has property $\rho$. Thus, if $\CM$ is an
%o-minimal expansion of a real closed field $R$ and $K$ is the
%algebraic closure of $R$ then it follows from Theorem
%\ref{interpret2} that every $\CM$-definable finite extension $\tilde
%G$ of a semisimple $K$-algebraic group $G$ is itself $K$-algebraic.
% Indeed, by the theorem, $\tilde G$ can be interpreted in $\la
%G,\cdot\ra$, after naming finitely many constants.

\section{Groups which are not definably connected}

In this section we prove an analogue of Theorem \ref{thm4} for definably compact
groups which are not assumed to be definably connected.

We assume that $\CM$ is a sufficiently saturated
o-minimal structure expanding a real closed field.
We still use $\widehat G$ to denote $G/G^{00}$.
Here are some preliminaries:

\begin{claim}\vlabel{notconnected} (i) Let $G$ be any group, and $H$ a subgroup of finite index. Then
there are $g_{1},..,g_{n}\in G$, and $h_{1},..,h_{m}\in H$ such that the structures
$(G, \cdot, H, g_{1},..,g_{n})$ and $(H,\cdot, a_{1},..,a_{n},h_{1},..,h_{m})$ are bi-interpretable, where
$a_{i}$ is conjugation by $g_{i}$.
\newline
(ii) In the special case when $G$ is definable in an $o$-minimal structure and $H = G^{0}$, then $G^{0}$
is definable in $\la G, \cdot \ra$ so in fact $(G,\cdot,g_{i})_{i}$ and $(G^{0},\cdot, h_{j},a_{i})_{j,i}$
are bi-interpretable.
\end{claim}
\proof
%(i) Since $G^0$ is divisible, every coset of $G^0$ in $G$
%contains a torsion element. Let $f_1,\ldots, f_k$ be a set of
%torsion elements representing these cosets and let $F\sub G$ be the
%finite group generated by these elements. Each $f\in F$ acts on
%$G^0$ by conjugation. The map $(g,f)\mapsto gf$ sends the semidirect
%product $H=G^0\rtimes F$ onto the group $G$. Its finite kernel is
%$\{(g,g^{-1}):g\in F\cap G^0\}$.

%For the interpretation in the opposite direction, note that
%$G^0=G^k$ (in the notation of \ref{div1}) for some $k$ hence $G^0$
%is interpretable in $G$. The bi-interpretability (over parameters)
%is easy to verify.
(i) is straightforward: let the $g_{i}$ be representatives of cosets of $H$ in $G$, and for each $i,k$ let
$g_{i}g_{k} = c_{ik}g_{r}$ for suitable $r$ and $c_{ik}\in H$. Let the $h_{j}$'s be the $c_{ik}$'s.
Details are left to the reader.
\newline
(ii) In
the notation of lemma \ref{div1}, there is an $n$ such that
$\sigma_n(G)\sub G^0$ and therefore, by the same lemma, there is a
$k$ such that $G^0=\sigma_n(G)\cdots \sigma_n(G)$ ($k$-times). This
implies that $G^0$ is definable (without parameters) in $G$.
\qed

\vspace{5mm}
\noindent
By the above, in order to understand an arbitrary definable group
$G$ we need to understand $G^0$ together with finitely many
definable automorphisms.

By \ref{structure}, every definably compact, definably connected
group $G$ is the almost direct product of the semisimple group
$[G,G]$ and $Z(G)^0$.

Clearly, every definable automorphism of $G$ leaves invariant both
$[G,G]$ and $Z(G)^0$, so we need to understand each of the two
groups, together with finitely many definable automorphisms.

Theorem \ref{theorem-ab} allows us to treat definable automorphisms
of a definable abelian group $A$ (by viewing their graph as a
subgroup of $A\times A$). Hence, we now need to examine  definable
automorphisms of definable semisimple groups.

\begin{claim} \vlabel{automorss} Let $\CM$ be an o-minimal expansion of a real closed field
$R$. If $G$ is an $R$-semialgebraic, definably connected, definably
compact semisimple group,  then every definable automorphism of $G$
is $R$-semialgebraic.
\end{claim}
\proof Assume first that $G$ is definably simple and $f:G\to G$ is a
definable automorphism. Because $G$ is definably compact it is
bi-interpretable with a real closed field $R_1$. Therefore, by
\cite{pps2}*{Proposition 4.8}, (a Borel-Tits-style result),
$f=g\circ h$, where $g$ is an $R$-semialgebraic automorphism of $G$
and $h$ is induced by an automorphism $\sigma$ of the semialgebraic field
$R_1$. The proof of Proposition 4.8 cited above shows that $\sigma$ is definable, if $f$ is definable. So
$\sigma$ is definable and thus the identity. So $f$ is semialgebraic, proving the claim in
the special case.

\vspace{2mm}
\noindent
Assume now that $G$ is semisimple and centerless. Hence, by
\ref{facts}, $G$ is definably isomorphic in $R$ to $H_1\times \cdots
\times H_n$, where each $H_i$ is a linear semialgebraic group,
defined over the real algebraic numbers $R_{alg}\sub R$. Without
loss of generality, $G=H_1\times \cdots \times H_n$, and we consider
each $H_i$ as a subgroup of $G$. For each $i=1,\ldots,k$ we denote
by $\pi_i:G\to H_i$  the projection map. Without loss of generality,
$dim H_n\leq \dim H_i$ for every $i\leq n$.

Let $f:G\to G$ be a definable automorphism. It is clearly sufficient
to see that $f:H_i\to G$ is semialgebraic for every $i=1,\ldots,
n$. We  prove that every definable embedding $\phi:H_i\to G$ is
semialgebraic. We first claim that $\pi_n\phi(H_i)=\{e\}$ or
$\pi_n\phi(H_i)=H_n$. Indeed, if $\pi_n\phi (H_i)\neq H_n$ then by
the dimension assumption, $ker \pi_n\phi \neq \{e\}$, which by
simplicity implies that $\pi_n\phi(H_i)=\{e\}$).

 Now, if $\pi_n\phi(H_i)=\{e\}$, we can reduce the problem by
induction to a group $G$ of smaller dimension. If
$\pi_n\phi(H_i)=H_n$ then the map $\pi_n\phi$ is an isomorphism of
$H_i$ and $H_n$ and therefore, similarly to the proof for the
definably simple case, $\phi|H_i$ is semialgebraic.

Let $G$ be an arbitrary definably connected semisimple group and
let $H=G/Z(G)$.  By \ref{T1.3}, $G$ is bi-interpretable with
$G/Z(G)$ over parameters in $G/Z(G)$ and $Z(G)$. Moreover, by the
last clause of Theorem \ref{interpret}, these parameters are
$0$-definable in $\la G,\cdot\ra$. Hence, there is a
definable group $G_1$ in the structure $\la G/Z(G), Z(G)\ra$, and a
$G$-definable isomorphism $\sigma:G\to G_1$, where $G_1$ and
$\sigma$ are invariant under every automorphism of the group $G$.

Now, let $f:G\to G$ be a definable automorphism in $\CM$. The map
$f$ induces an automorphism $f_1$ of $G/Z(G)$ and since $G/Z(G)$ is
centerless ($Z(G)$ is finite) it follows from the centerless case
that $f_1$ is $R$-semialgebraic. As we just pointed out, $f$ leaves
$G_1$ invariant and hence $f_1|G_1$ is an $R$-semialgebraic
automorphism of $G_1$. Composing with the $f$-invariant $\sigma$ we
see that $f$ itself is also $R$-semialgebraic.\qed

%{\bf The above was quite complicated and subtle. Is this the best we
%can do?}
%  It's OK. AP

The following lemma is general.

\begin{lemma}\vlabel{automquot} Let $\bf G$ be an arbitrary group,
 $A\sub G$  a central subgroup of $G$, such that for some number
 $k$, every element of
$G/A$ equals the product of $k$ commutators from $G/A$.

Let $f:G\to G$ be a group automorphism of $G$ such that
$f(A)=A$. Then $f$ is definable in the structure ${\bf G}=\la
G,,\cdot,A,f|A, f|G/A\ra$.
\end{lemma}
\proof We use Beth definability theorem: Namely, we take $\la \tilde G,
\tilde A, \tilde f|A, \tilde f|G/A\ra$ elementarily equivalent to
$\bf G$ and show that there is unique automorphism $g:\tilde G\to
\tilde G$ leaving $A$ invariant such that $g|A=\tilde f|A$ and
$g|G/A=\tilde f|G/A$. Note that the assumption on $G$ implies that
$\tilde G$ is still the product of $\tilde A$ and $[\tilde G,\tilde G]$
(this is true in $G$ and because $[G/A,G/A]$ is
generated in finitely many steps, it becomes a first order statement
true in $\tilde G$ as well).

Assume that we have $g,h:\tilde G\to \tilde G$ automorphisms as above
and consider $gh^{-1}$. Then $gh^{-1}|\tilde A=id$ and
$gh^{-1}|\tilde G/\tilde A=id$. We may therefore assume that
$g|\tilde A$ and $g|\tilde G/\tilde A$ are the identity maps and aim
to show that $g=id$.

Because $g|\tilde G/\tilde A=id$, for every $x\in \tilde G$, we have
$x^{-1}g(x)\in \tilde A$ and hence there exists a function $a:\tilde
G\to \tilde A$ such that $g(x)=xa(x)$. We claim that $a$ is a group
homomorphism: For $x,y\in \tilde G$ we have
$$xya(xy)=g(xy)=g(x)g(y)=xa(x)ya(y)=xya(x)a(y)$$(because $a(x),
a(y)$ are central elements). It follows that $a(xy)=a(x)a(y)$.

For every $x\in A$ we have $g(x)=x$, hence $a(x)=e$. Also, if
$b=xyx^{-1}y^{-1}$ is a commutator in $\tilde G$ then
$g(b)=ba(xyx^{-1}y^{-1})=b$, hence $a(b)=e$. But then $ker(a)$
contains both $A$ and the commutator subgroup of $\tilde G$. Because
$\tilde G$ is generated by these two groups, $a(x)=e$ for all $x\in
\tilde G$ and therefore $g=id$.\qed

%{\bf Was the use of Beth really necessary? It seems like there
%should be an explicit definition for $f$}
% It's OK. AP

\begin{theorem} If $G$ is a definably compact group in an o-minimal expansion
$\CM$ of an ordered group then it is elementarily equivalent to a
definably compact, semialgebraic (over parameters) group $H$ over a real closed field,
with $\dim H=\dim G$.
\end{theorem}
\proof As we saw above, $G$
 is bi-interpretable, over parameters
from $G^0$, with $G^0$ together with the action of finitely many
automorphisms $f_1,\ldots, f_k$. For simplicity, we denote $G^0$ by
$H$ and $Z(H)^0$ by $A$.

We also saw, in \ref{automquot}, that the structure $$\la H, \cdot,
\{f_1, \ldots,f_k\}, \{c_1,\ldots, c_r\}\ra $$ (for constants
$c_1,\ldots, c_r\in H$) is definable in $$\la H,\cdot, A,
\{f_i|A:i=1,\ldots, k\},\{f_i|H/A:i=1,\ldots, k\},
\{c_1,\ldots, c_r\}\ra.$$

We denote each $f_i|H/A$ by $g_i$ and each $f_i|A$ by $h_i$.

 By
Theorem \ref{thm1}, the group $H$  is definable, over
parameters, in the two-sorted structure $\la H/A,  A \ra$. Putting it
all together, we see that $G$ is definable, over parameters, in
$\la H/A,\{g_i:i=1,\ldots, k\}, A, \{h_i:i=1,\ldots, k\}\ra$. (where
$H/A$ and $A$ are endowed with their group structure).

By \ref{T1.3} the semisimple group $H/A$ is definably isomorphic to
a semialgebraic group $G_0$ over $R_{alg}\sub R$, for a real closed
field $R$ and by \ref{automorss}, each $g_i$ is sent by this
isomorphism to an $R$-semialgebraic automorphism of $G_0$, possibly
defined over parameters.

 The structure $\la A,+, \{h_1,\ldots,h_k\}\ra $ is clearly a
 reduct of the structure ${\bf A}_{ab}$ considered in Theorem \ref{theorem-ab}
 (since every automorphism of $G$ gives rise to a subgroup of $G\times G$).
 Therefore, it is elementarily equivalent to
 an expansion of a connected, compact, abelian real Lie
 group $\widehat A$ (with $\dim \widehat A=\dim A$),
  by Lie group automorphisms $\widehat h_1,\ldots,
 \widehat h_k$. Finally, $\widehat A$, as a compact Lie group,  is
 isomorphic to a real algebraic linear group $L$. This isomorphism takes
 each $\widehat h_i$ to a a Lie subgroup of $L^2$, which itself must
 be semialgebraic (indeed, this last fact follows for example from
 \cite{pps1}*{3.3}, applied to the o-minimal structure
 $\reals_{an}$, in which every definable compact linear group is
 definable).

Hence, by going to a sufficiently saturated real closed field
$\widetilde R$, we can find constants $d_1,\ldots, d_k\in \widetilde
R$  such that $$\CM_1=\la A,\{h_i:i=1,\ldots,k\},
H/A,\{g_i:i=1,\ldots,k\},
 \{c_i:=i=1,\ldots,r\}\ra$$ is elementarily
equivalent to
$$\CM_2=\la L(\widetilde
R),\{\widehat h_i:i=1,\ldots,k\},  G_0(\widetilde R),\{\widehat
g_i:i=1\ldots,k\}, \{d_i:i=1,\ldots,r\}\ra,$$ with $G_0$, $L$, and
the automorphisms $\widehat g_i, \widehat h_i$ all semialgebraic.

Because $G$ is definable over parameters  in $\CM_1$, it is
elementarily equivalent to a  group definable (over parameters)
in $\CM_2$, and this last group must be semialgebraic.\qed
\\

\noindent{\bf Remark.} By Lemma \ref{abelianclaim}, the parameters
in $A$ can be realized as real algebraic elements in an elementarily
equivalent real group. However, we don't know how to do the same for
the parameters in $H/A$.

\section{Compact domination for definably compact groups}
Here we give another application of Corollary 6.4. The ``compact domination conjecture" for
{\em definably compact}, definably connected groups in (saturated) $o$-minimal expansions of real closed
fields, was introduced in \cite{nip1}.
The conjecture says that $G$ is dominated by $G/G^{00}$ equipped with its Haar measure.
Namely, writing $\pi:G\to G/G^{00}$ for the canonical surjective homomorphism, for any definable subset $X$
of $G$ the set of $c\in G/G^{00}$ such that $\pi^{-1}(c)$ intersects both $X$ and its complement, has Haar
measure $0$.
We sometimes just say ``$G$ is compactly dominated".
The conjecture was proved in \cite{nip1} for $G$ with ``very good reduction", and by part (ii) of Theorem
\ref{T1.3} of the current paper, this is the case for semisimple definably connected groups. In \cite{nip2}
compact domination was proved for $G$ commutative.
With \ref{structure} we know that arbitrary $G$ (definably compact, definably connected) almost splits into
its semisimple and abelian parts, and one would expect that this makes it easy  to deduce compact
domination of $G$ from the two special cases.

We first prove:
\begin{theorem} Every definably connected, definably compact group in $\CM$ a sufficiently saturated
expansion of a real closed
field is compactly dominated.
\end{theorem}
\pf We first prove the result for a group $G\times H$, with $G$ commutative and $H$ semisimple.
  It is sufficient to prove the result for each definable set separately so we may assume
that the language is countable. We fix $\CM_0\sub \CM $ a countable model. We let
$$\pi:G\times H\to G/G^{00}\times H/H^{00},$$ with $\pi=(\pi_1,\pi_2)$ and
$\pi_1:G\to G/G^{00}$, $\pi_2:H\to H/H^{00}$.

The following observation is true in greater generality (for any type-definable equivalence relation),
but we only observe it in the o-minimal setting: If $K$ is a definably compact group in $\CM$, definable
over $M_0$  and $a_1,a_2\in K$
realize the same type over $M_0$ then the lie in the same $K^{00}$-coset.

Indeed, if $\sigma$ is any automorphism of $\CM$ then it induces a continuous
(with respect to the logic topology) automorphism of $K/K^{00}$ which fixes all the torsion points of $K$
(since they belong to $M_0$). But $\pi(Tor(K))$ is dense in $K/K^{00}$, therefore $\sigma$ induces the
identity map on $K/K^{00}$. This is clearly sufficient.

Let $X\sub G\times H$ be a definable set over $M_0$ and assume, towards contradiction, that
the set $$B=\{(g',h')\in G/G^{00}\times H/H^{00}: \pi^{-1}(g',h')\cap X\neq \emptyset \,\&\,
\pi^{-1}(g',h')\cap X^c\neq \emptyset\}$$ has positive measure.

 Recall from \cite{nip2} that ${\CM}^{*}$ denotes the expansion of $\CM$ by adjoining relations for all
 externally definable sets, and that the theory of this expansion is weakly o-minimal. From \cite{nip1},
 $Fin$ denotes the finite elements of $\CM$ and $Inf$ the infinitesimals, both definable in ${\CM}^{*}$.
 $Fin/Inf$ identifies with $\reals$ and the structure on it that is induced from ${\CM}^{*}$ is o-minimal.
 Moreover $H/H^{00}$ is a definable subset of some $(Fin/Inf)^{n}$ and the logic topology on $H/H^{00}$
 coincides with its topology as a definable group in $Fin/Inf$. In particular any subset of $H/H^{00}$
 definable in ${\CM}^{*}$ which has positive Haar measure, is of maximal o-minimal dimension, so has
 interior.
The results in \cite{nip2} give a similar picture for $G/G^{00}$. Namely $G^{00}$ is definable in
${\CM}^{*}$ and
$G/G^{00}$ is semi-o-minimal, namely lives in the product of finitely many definable o-minimal sets in
${\CM}^{*}$. The two topologies (logic, o-minimal) again coincide, so definable (in ${\CM}^{*}$) sets of
positive Haar measure have interior.

The group $G/G^{00}\times H/H^{00}$, with all its induced $\CM^*$ structure,  is also
semi-o-minimal and because the set $B$ is definable in $\CM^*$ there are open sets $U\sub G/G^{00}$ and
$V\sub H/H^{00}$
with $U\times V \sub B$. We  claim that there exists $g'\in U$ such that all elements of $\pi_1^{-1}(g')$
realize the same type in $\CM$, over $M_0$. Indeed, because $G$ is compactly dominated, for every
$M_0$-definable subset $X$ of $G$, the set of all $g'\in G/G^{00}$ such that $\pi_1^{-1}(g')\cap X\neq
\emptyset$ and $\pi_1^{-1}(g')\cap X^c\neq \emptyset$ has Haar measure zero. So, after removing countably
many such sets (each of measure zero), the pre-image of every $g'\in U$ under $\pi_1$ is contained in a
complete $\CM$-type over $M_0$.

We fix one such $g'\in U$ as above, $g\in \pi_1^{-1}(g')$, and consider the set
$X_g=\{h\in H:(g,h)\in X\}$. We claim that for every $h'\in V$, the sets $\pi_2^{-1}(h')\cap X_g$ and
$\pi_2^{-1}(h')\cap X_g^c$ are both nonempty, contradicting the fact that $H$ is compactly dominated.
Indeed, if $h'\in V$ then, by assumption on $B$, there are $g_1,g_2\in \pi_1^{-1}(g')$ and $h_1,h_2\in
\pi_2^{-1}(h')$ with $(g_1,h_1)\in X$ and $(g_2,h_2)\in X^c$. Because $g_1,g_2$ and $g$ all realize
the same type over $M_0$, there are $h_3,h_4\in H$ conjugates over $M_0$ of $h_1,h_2$, respectively, with
$(g,h_3)\in X$ and $(g,h_4)\in X^c$. By our earlier observation, $h_3$ and $h_4$ belong to the pre-image of
$h'$, so $\pi_2^{-1}(h')\cap X_g$ and $\pi_2^{-1}(h')\cap X_g^c$ are non-empty. Contradiction.
 We thus showed that $G\times H$ is compactly dominated.

 The result for an arbitrary definably compact group follows from the special case
 using Corollary \ref{structure}, noting that compact domination is preserved under quotients
 (using the fact that a definable surjective homomorphism $\sigma:G_1\to G_2$ of definably
 compact groups sends $G_1^{00}$ onto $G_2^{00}$).

\qed

\section{Appendix: On Abelian groups}

Since all groups here are abelian we write them additively. The
first two lemmas are standard and we include them for
completeness.

\begin{lemma}\vlabel{abelian1} Let $B\sub A$ be two abelian
divisible groups such that $B$ contains all torsion elements of $A$
and let $A_0$ be a subgroup of $A$. Then the following are
equivalent:
\\(i) $A_0$ is a maximal subgroup of $A$ such
that $A_0\cap B=\{0\}$.
\\(ii) $A=A_0\oplus B$ and $A_0$ is divisible.
\\(iii) $A=A_0\oplus B$.
\end{lemma}
\proof We only need to see that (i) implies (ii).
 Notice that the assumptions on $B$ (with $A$ divisible) are equivalent to:
 for every $n\geq 1$ and $a\in A$, if $na\in B$ then
$a\in B$.

 We first show:

\noindent{\bf Claim } For every $n\geq 1$ and $a\in A$, if $na\in
A_0+B$ then $a\in A_0+B$.
\\

We use induction on $n$: The case $n=1$ is obvious. For general $n$,
if $na=a_0+b$ for some $a_0,\in A_0, b\in B$, then, since $A$ is
divisible there exists $a'\in A$ such that $na'=a_0$. It follows
that $n(a-a')\in B$ and therefore, by the above observation,
$a-a'\in B$.  It is therefore sufficient to prove that $a'\in
A_0+B$.

If $a'\notin A_0$ then, by the maximality of $A_0$, there exists
$k\in \mathbb Z$, and $a_1\in A_0$ such that $ka'+a_1=b'\in B$ and
$b' \neq 0$. We write $k=mn+\ell$, $0\leq \ell<n$ and then we have
$mna'+\ell a'+a_1=b'$. If $\ell=0$ then $mna'+a_1=b'\neq 0$, which
is impossible because $mna'+a_1\in A_0$. It follows that $\ell\neq
0$ and we have $\ell a'=(-mna'-a_1)+b\in A_0+B$. By induction, we
have $a'\in A_0+B$. End of Claim.
\\

We can now prove (ii): For $a\in A$, if $a\notin A_0$ then, by the
maximality of $A_0$, there exist $k\in \mathbb Z$ and $a_0\in A_0$
such that $0\neq ka+a_0=b\in B$. By assumption on $A_0$, we have
$k\neq 0$ (and  $ka\in A_0+B$). The last claim implies that $a\in
A_0+B$.

To see that $A_0$ is divisible, take $a_0\in A_0$ and $n\geq 1$.
Because $A$ is divisible, we have  $a_1\in A_0, b\in B$, such that
$na_1+nb=a_0$. It follows that $nb=0$ and then we also have
$na_1=a_0$. \qed

\begin{lemma} \vlabel{abelian2} Let $A,B$ be two divisible abelian
groups such that $B$ has unbounded exponent. Assume that $\phi:B\to
A$ is a group embedding, with $Tor(A)\sub \phi(B)$. Then $\phi$ is
an elementary map.
\end{lemma}
\proof We assume that $B\subsetneq A$ and by moving to an elementary
extension of $\la A,+,B\ra$ we may assume that $B$ contains a
torsion-free element. By \ref{abelian1}, we can write
$B=Tor(A)\oplus B_0$ for some torsion-free divisible subgroup $B_0$.
Because $B$ has unbounded exponent, $B_0\neq \{0\}$. Applying again
\ref{abelian1}, we can write $A=Tor(A)\oplus B_0\oplus A_0$ for
$A_0\neq \{0\}$. We now use the fact that $B_0\prec(B_0\oplus A_0)$
(as divisible torsion-free abelian groups) to conclude that $B\prec
A$.\qed

We also need the following claim on definable abelian groups in
o-minimal structures.

 \begin{lemma}\vlabel{abelianclaim} Let $\CM$ be an o-minimal expansion of an ordered
 group. If
$G$ is a definable abelian group (possibly not definably connected)
and $C$ is a finite subset of $G$ then $\la G,+, \{b:\in C\}\ra $
(namely, we add a constant to every element of $C$) is elementarily
equivalent to a real algebraic group of the same dimension, with
finitely many real algebraic elements named.
\end{lemma}

\proof  Assume first that  $G$ is definably connected. It follows
that it is divisible.

By \cite{ps}, there exists $G_0\sub G$ a torsion-free definable
subgroup of $G$ with $G/G_0$ definably compact.  Because $G_0$ is divisible and torsion-free it is
elementarily equivalent to $\reals^{\dim G_0}$. The group
$G_1=G/G_0$ is a definably compact, definably connected group and
therefore by \cite{edot} (and, in the case that $\CM$ expands an
ordered group also by \cite{es} and \cite{pe}), $Tor(G_1)$ is
isomorphic to the torsion group of the real torus ${\mathbb T}^{\dim
G_1}$. It follows (say, by \ref{abelian2}) that $G_1$ is
elementarily equivalent to the semialgebraic ${\mathbb T}^{\dim
G_1}$ and $G$ is elementarily equivalent to the group $\reals^{\dim
G_0}\times {\mathbb T}^{\dim G_1}$.

Finally, consider the divisible hull $H$ in $G$ of the group
generated by the set $C$ and $Tor(G)$. By \ref{abelian1}, $H$ can be
written as the direct sum of $Tor(G)$ and $\mathbb
Qc_1\oplus\cdots\oplus \mathbb Q c_k$, for some $c_1,\ldots,c_k\in
C$. By \ref{abelian2}, $H$ is an elementary subgroup of $G$. It can
be realized also as an elementary subgroup of $\reals^{\dim
G_0}\times {\mathbb T}^{\dim G_1}$. Moreover, since all torsion
elements are real algebraic and we can also choose real algebraic
elements which are torsion-free and $\mathbb Q$-independent (using the well-known fact that
the field of real algebraic numbers is infinite-dimensional as a $\mathbb Q$-vector space), it
follows that
 $\la G,+,\{g:g\in C\cup Tor(G)\}\ra$ is elementarily
equivalent to $\reals^{\dim G_0}\times {\mathbb T}^{\dim G_1}$, with
names for all torsion elements and finitely many other elements, all
real algebraic.

Assume now that  $G$ is not definably connected.  Then it equals a
direct sum of its connected component $G^0$ and a finite group and
therefore, by the above, it is elementarily equivalent to a
semialgebraic group $H$ of the same dimension, which can be defined
over the real algebraic numbers. We can handle similarly finitely
many named elements in $G$.\qed
\\

 Given an expansion $\bf
A$ of an abelian group $A$, consider a sub-language $L_{ab}$ which
has a predicate $R_S$ for every $0$-definable (in ${\bf A}$)
subgroup $S\sub A^n$, $n\in \mathbb N$, as well as symbols for $+$ and $0$. Let ${\bf A}_{ab}$ denote the
reduct of $\bf A$ to $L_{ab}$. Such a structure has been sometimes called
an {\em abelian structure}. We call the subgroups $S$ above the {\em basic} ones in the structure
${\bf A}_{ab}$. If $B$ is a subgroup of $A$ then we
denote by ${\bf B}_{ab}$ the $L_{ab}$-induced structure on $B$, namely the interpretation of $R_{S}$ is
just its intersection with $B^{n}$. The next fact is a restatement of \ref{module1}, but we add a bit more
information in item 1.

\begin{fact}\vlabel{module1.1} In the above setting (no o-minimality is assumed)
\begin{enumerate}
\item (i) The theory of the structure ${\bf A}_{ab}$ eliminates quantifiers (in the language
$L_{ab}$).
\newline
(ii)  Moreover $Th({\bf A}_{ab})$ is axiomatized as follows:
(a) Axioms for abelian groups, (b) Each symbol $R_{S}$ denotes a subgroup, (c) Axioms for the defining
properties of $R_S$:
 If $S_{1}$ is a projection of $S_{2}$ then  $R_{S_{1}}$ denotes the corresponding
projection of $R_{S_{2}}$, and if $S_1=\{x\in A^n: S_2(x,0)\}$ then $R_{S_1}$ denotes
 the corresponding fiber of $R_{S_2}$, (d) Axioms about the index (a given finite number or $\infty$)
 of $R_{S_{1}}$ in $R_{S_{2}}$ whenever $S_{1}\leq S_{2} \leq A^{n}$ are basic.

\item Assume that $B\leq A$ is a subgroup of $A$.

Then ${\bf B}_{ab}\prec {\bf A}_{ab}$ if and only if the following
hold: \newline
(i) For every $0$-definable (in ${\bf A}_{ab}$)
subgroup $S\leq A^{n+k}$ and $b\in B^k$,
$$S(B^n,b)\neq \emptyset \Leftrightarrow S(A^n,b)\neq
\emptyset.$$ (ii) For all $0$-definable (in ${\bf
A}_{ab}$) subgroups $S_1 \leq S_2 \leq A^n$,
$$[S_2:S_1]= [S_2\cap B^n:S_1\cap B^n],$$ with the meaning that if
this index is infinite on one side  then it is infinite on the
other.
\item Assume that ${\bf A}_{ab}$ has DCC on $0$-definable
subgroups. Then, for every ${\bf B}_{ab}\prec{\bf A}_{ab}$ there
exists a surjective group homomorphism $\phi:A\to B$ which is the
identity map on $B$ and in addition sends every $0$-definable $S\sub
A^n$ onto $S\cap B^n$. (We call such a $\phi$ {\em a homomorphic
retract}).
\end{enumerate}
\end{fact}
\proof (1)(i) is proved in \cite{GR}.
\newline
1(ii) can be extracted from the proof of the quantifier elimination result in
\cite{GR}, in exactly the same way as the analogous statement for theories of modules is deduced from the
proof of $pp$ elimination for modules. See Theorem 1.1 and Corollary 1.5 of \cite{Z}. In fact in the
statements on indices only basic subgroups of $A$ itself (rather than $A^{n}$) need be considered.

(2) follows from (1).

(3) Using the quantifier elimination result above, the proof of
(3) is basically identical to that of Theorem 2.8, p.28, in
\cite{Prest}.\qed
\\


\begin{bibdiv}
\begin{biblist}
\normalsize


\bib{b1}{article}{
   author={Berarducci, Alessandro},
   title={Zero groups and maximal tori},
   journal={Logic Colloquium 04, ed. A. Andretta, K. Kearnes and  D. Zambella},
   volume={Lecture notes in Logic 29},
   %publisher={},
   %place={},
   date={2006},
   pages={33-47},
%   isbn={0-521-59838-9},
%   review={\MR{1633348 (99j:03001)}},
}




\bib{b2}{article}{
   author={Berarducci, Alessandro},
   title={O-minimal spectra, infinitesimal subgroups and cohomology},
   journal={J. Symbolic Logic},
   volume={72},
   %publisher={},
   %place={},
   date={2007},
   number={4},
   pages={1177-1193},
%   isbn={0-521-59838-9},
%   review={\MR{1633348 (99j:03001)}},
}

\bib{bo}{article}{
author = {Berarducci, Alessandro},
author =  {Otero, Margarita},
title = {O-minimal fundamental groups, homology and manifolds},
journal = {J. London Math. Soc.},
date = {2002},
number = {65},
pages = {257-270},
}


\bib{vdd}{book}{
   author={van den Dries, Lou},
   title={Tame topology and o-minimal structures},
   series={London Mathematical Society Lecture Note Series},
   volume={248},
   publisher={Cambridge University Press},
   place={Cambridge},
   date={1998},
   pages={x+180},
%   isbn={0-521-59838-9},
%   review={\MR{1633348 (99j:03001)}},
}



\bib{e1}{article}{
   author={Edmundo, Mario J.},
   title={Solvable groups definable in o-minimal structures},
   journal={J. Pure and Appl. Algebra},
   volume={185},
   %publisher={Cambridge University Press},
   %place={Cambridge},
   number={1-3},
   date={2003},
   pages={103\ndash 45},
%   isbn={0-521-59838-9},
%   review={\MR{1633348 (99j:03001)}},
}

\bib{e2}{article}{
   author={Edmundo, Mario J.},
   title={Covers of groups definable in o-minimal structures},
   journal={Illinois J. Math.},
   volume={49},
      number={1},
   date={2005},
   pages={99\ndash 120},
%   isbn={0-521-59838-9},
%   review={\MR{1633348 (99j:03001)}},
}

\bib{e3}{article}{
   author={Edmundo, Mario J.},
   author={Pantelis Eleftheriou}
   title={The universal covering homomorphism in o-minimal expansions of groups},
   journal={Math. Logic Quart.},
   volume={53},
      number={6},
   date={2007},
   pages={571\ndash 582},
%   isbn={0-521-59838-9},
%   review={\MR{1633348 (99j:03001)}},
}

\bib{ejp}{article}{
   author={Edmundo, M{\'a}rio J.},
   author={Jones, J},
   author={Peatfield, N.}
   title={Hurewicz Theorems for definable groups, Lie groups and their cohomologies},
   journal={preprint},
   volume={},
   date={},
   number={},
   pages={},
%   issn={0219-0613},
%   review={\MR{2114966 (2005m:03073)}},
}




\bib{edot}{article}{
   author={Edmundo, M{\'a}rio J.},
   author={Otero, Margarita},
   title={Definably compact abelian groups},
   journal={J. Math. Log.},
   volume={4},
   date={2004},
   number={2},
   pages={163--180},
%   issn={0219-0613},
%   review={\MR{2114966 (2005m:03073)}},
}

\bib{es}{article}{
   author={Eleftheriou, Pantelis E.},
   author={Starchenko, Sergei},
   title={Groups definable in ordered vector spaces over ordered division
   rings},
   journal={J. Symbolic Logic},
   volume={72},
   date={2007},
   number={4},
   pages={1108--1140},
%   issn={0022-4812},
%   review={\MR{2371195 (2008j:03051)}},
}


\bib{GR}{article}{
   author={Gute B. Hans},
   author={Reuter K. K},
     title={The last word on elimination of quantifiers in modules},
   journal={J. Symbolic Logic},
   volume={55},
   date={June 1990},
   number={2},
   pages={670-673},
   %issn={0168-0072},
   %review={\MR{1729742 (2000m:03090)}},
}


    \bib{HM}{book}{
   author={Hofmann H. Karl},
   author={Morris A. Sidney},
   title={The structure of compact groups},
   series={de Gruyter Studies in Mathematics 25},
   %note={Translated from the Russian and with a preface by D. A. Leites},
   publisher={de Gruyter},
   place={Berlin},
   date={1998},
   pages={xvii+835},
  % isbn={3-540-50614-4},
   % review={\MR{1064110 (91g:22001)}},
}

\bib{nip1}{article}{
   author={Hrushovski, Ehud},
   author={Peterzil, Ya'acov},
   author={Pillay, Anand},
   title={Groups, measures, and the NIP},
   journal={J. Amer. Math. Soc.},
   volume={21},
   date={2008},
   number={2},
   pages={563--596},
%   issn={0894-0347},
 %  review={\MR{2373360 (2008k:03078)}},
}

\bib{nip2}{article}{
 author={Hrushovski, E.},
      author={Pillay, A.},
     title={NIP and invariant measure},
   journal={preprint},
   volume={},
   date={2007},
   number={},
   pages={},
%   issn={0002-9947},
 %  review={\MR{1707202 (2001b:03036)}},
}


\bib{opp}{article}{
 author={Otero, M.},
 author={Peterzil, Y.}
      author={Pillay, A.},
     title={Groups and rings definable in o-minimal expansions of real closed fields},
   journal={Bull. LMS},
   volume={28},
   date={1996},
   number={},
   pages={7-14},
%   issn={0002-9947},
 %  review={\MR{1707202 (2001b:03036)}},
}



    \bib{OV}{book}{
   author={Onishchik, A. L.},
   author={Vinberg, {\`E}. B.},
   title={Lie groups and algebraic groups},
   series={Springer Series in Soviet Mathematics},
   note={Translated from the Russian and with a preface by D. A. Leites},
   publisher={Springer-Verlag},
   place={Berlin},
   date={1990},
   pages={xx+328},
  % isbn={3-540-50614-4},
   % review={\MR{1064110 (91g:22001)}},
}

\bib{pe}{article}{
   author={Peterzil, Y.},
     title={Returning to semi-bounded sets},
   journal={to appear in JSL},
   volume={},
   date={2007},
   number={},
   pages={},
   %issn={0168-0072},
   %review={\MR{1729742 (2000m:03090)}},
}
        \bib{pps1}{article}{
   author={Peterzil, Y.},
   author={Pillay, A.},
   author={Starchenko, S.},
   title={Definably simple groups in o-minimal structures},
   journal={Trans. Amer. Math. Soc.},
   volume={352},
   date={2000},
   number={10},
   pages={4397--4419},
%   issn={0002-9947},
 %  review={\MR{1707202 (2001b:03036)}},
}

    \bib{pps2}{article}{
   author={Peterzil, Y.},
   author={Pillay, A.},
   author={Starchenko, S.},
   title={Simple algebraic and semialgebraic groups over real closed fields},
   journal={Trans. Amer. Math. Soc.},
   volume={352},
   date={2000},
   number={10},
   pages={4421--4450 (electronic)},
 %  issn={0002-9947},
 %  review={\MR{1779482 (2001i:03080)}},
}

\bib{pps3}{article}{
   author={Peterzil, Y.},
   author={Pillay, A.},
   author={Starchenko, S.},
   title={Linear groups definable in o-minimal structures},
   journal={J. Algebra},
   volume={247},
   date={2002},
   number={1},
   pages={1--23},
  % issn={0021-8693},
  % review={\MR{1873380 (2002i:03043)}},
}

\bib{ps1}{article}{
   author={Peterzil, Ya'acov},
   author={Starchenko, Sergei},
   title={Definable homomorphisms of abelian groups in o-minimal structures},
   journal={Ann. Pure Appl. Logic},
   volume={101},
   date={2000},
   number={1},
   pages={1--27},
   %issn={0168-0072},
   %review={\MR{1729742 (2000m:03090)}},
}
    \bib{ps}{article}{
   author={Peterzil, Ya'acov},
   author={Steinhorn, Charles},
   title={Definable compactness and definable subgroups of o-minimal groups},
   journal={J. London Math. Soc. (2)},
   volume={59},
   date={1999},
   number={3},
   pages={769--786},
%   issn={0024-6107},
 %  review={\MR{1709079 (2000i:03055)}},
}

\bib{pi}{article}{
   author={Pillay, A},
   title={On groups and fields definable in o-minimal structures},
   journal={J. Pure Appl. Algebra},
   volume={53},
   date={1988},
   number={3},
   pages={239--255},
%   issn={0024-6107},
 %  review={\MR{1709079 (2000i:03055)}},
}


\bib{Prest}{book}{
   author={Prest, Mike},
   title={Model Theory and Modules},
   series={London Mathematical Society Lecture Note Series},
   volume={130},
   publisher={Cambridge University Press},
   place={Cambridge},
   date={1988},
   pages={x+180},
%   isbn={0-521-59838-9},
%   review={\MR{1633348 (99j:03001)}},
}

\bib{Z}{article}{
   author={Ziegler, Martin},
       title={Model theory of modules},
   journal={Ann. Pure Appl. Logic},
   volume={26},
   date={1984},
   number={},
   pages={149-213},
   %issn={0168-0072},
   %review={\MR{1729742 (2000m:03090)}},
}


\end{biblist}
\end{bibdiv}
\end{document}